\documentclass{amsart}
\usepackage{amsfonts}

\usepackage{graphicx}
\usepackage{amscd}
\usepackage{amsmath}

\newtheorem{theorem}{Theorem}
\theoremstyle{plain}

\newtheorem{definition}[theorem]{Definition}

\newtheorem{lemma}[theorem]{Lemma}

\newtheorem{fact}[theorem]{Fact}

\newtheorem{proposition}[theorem]{Proposition}
\theoremstyle{definition}
\newtheorem{remark}[theorem]{Remark}

\renewcommand\Bbb{\mathbb} 

\newcommand\set[1]{\ensuremath{\{#1\}}}
\newcommand\linebreakx{\linebreak}
\newenvironment{romenumerate}{\begin{enumerate}% gives (i), (ii) etc.
 }{\end{enumerate}}

\newcounter{thmenumerate}
\newenvironment{thmenumerate}
{\setcounter{thmenumerate}{0}%
 \def\item{\par% \ifnum\thethmenumerate=0\else\par\fi %\noindent\fi
 \refstepcounter{thmenumerate}\textup{(\roman{thmenumerate})\enspace}}
}
{}

\begin{document}
\title[Complex interpolation of compact operators]
{Complex interpolation of compact operators mapping into the couple $%
(FL^{\infty },FL_{1}^{\infty })$.}

\author[M. Cwikel]{Michael Cwikel}
\address{Department of Mathematics, Technion - Israel Institute of
  Technology, Haifa 32000, Israel}
\email{mcwikel@{math.technion.ac.il}}
\urladdr{http://www.math.technion.ac.il/\~{}mcwikel/}

\author[S. Janson]{Svante Janson}
\address{Department of Mathematics, Uppsala University, P.O.\ Box 480, S-751 06
Uppsala, Sweden}
\email{svante.janson{@math.uu.se}}
\urladdr{http://www.math.uu.se/\~{}svante/}

\thanks{The research of the first named author was supported by the
  Technion V.P.R.\ 
Fund and by the Fund for Promotion of Research at the Technion.}
\date{June 22, 2006}

\begin{abstract}

If $(A_{0},A_{1})$ and $(B_{0},B_{1})$ are Banach couples and a linear
operator $T:A_{0}+A_{1}\rightarrow B_{0}+B_{1}$ maps $A_{0}$ compactly
into $B_{0}$ and maps $A_{1}$ boundedly into $B_{1}$, does $T$ necessarily
also map $[A_{0},A_{1}]_{\theta }$ compactly into $[B_{0},B_{1}]_{\theta
}$ for $\theta \in (0,1)$?

After 42 years this question is still not answered, not even in the case
where $T:A_{1}\rightarrow B_{1}$ is also compact. But affirmative answers
are known for many special choices of $(A_{0},A_{1})$ and
$(B_{0},B_{1})$. Furthermore it is known that it would suffice to resolve
this question in the special case where $(B_{0},B_{1})$ is the special
couple $(\ell ^{\infty }(FL^{\infty }),\ell ^{\infty }(FL_{1}^{\infty
}))$. Here $FL^{\infty }$ is the space of all sequences $\{\lambda
_{n}\}_{n\in \Bbb{Z}}$ which are Fourier coefficients of essentially
bounded functions, and $FL_{1}^{\infty }$ is the weighted space of all
sequences $\{\lambda _{n}\}_{n\in \Bbb{Z}}$ such that $\{e^{n}\lambda
_{n}\}_{n\in \Bbb{Z}}\in FL^{\infty }$.

We provide an affirmative answer to this question in the related but
simpler case where 
$(B_{0},B_{1})$ is the special couple $(FL^{\infty },FL_{1}^{\infty })$.
\end{abstract}

\maketitle

\section{Introduction}

\smallskip This paper deals with the following question, which originated
with Alberto Calder\'{o}n's study \cite{Ca} of complex interpolation spaces,
and has been open for the 42 years which have elapsed since the writing of 
\cite{Ca}.

\textbf{\textit{Question C: }}\textit{Suppose that }$A_{0}$\textit{\ and }$%
A_{1}$\textit{\ are compatible Banach spaces, i.e., they form a Banach pair,
and that so are }$B_{0}$\textit{\ and }$B_{1}$\textit{. Suppose that }$%
T:A_{0}+A_{1}\rightarrow B_{0}+B_{1}$\textit{\ is a linear operator such
that }$T:A_{0}\rightarrow B_{0}$\textit{\ compactly and }$T:A_{1}\rightarrow
B_{1}$\textit{\ boundedly. Does it follow that }$T$\textit{\ maps the
complex interpolation space }$[A_{0},A_{1}]_{\theta }$\textit{\ into the
compact interpolation space }$[B_{0},B_{1}]_{\theta }$\textit{\ compactly
for each }$\theta \in (0,1)$\textit{?}

\smallskip

Affirmative answers have been obtained over the years for a considerable
number of particular cases of Question C. See, for example, \cite{Ca}, \cite
{LP}, \cite{P}, \cite{CKS}, \cite{CP}, \cite{Cw2}, and \cite{CwK}. More
recent studies of other aspects of this question can be found in \cite{CFM}, 
\cite{CKM}, \cite{CKM2} and \cite{Sc}. We also refer to \cite{CwJ} and the
website \cite{mcweb} for some general remarks and further details about this
question and related questions, including a remarkable ``almost
counterexample'' found by Fedor Nazarov.

\smallskip Although we still cannot answer Question C, we can, apparently
for the first time in nearly thirteen years, i.e., since the submission of 
\cite{CwK}, further enlarge the family of partial affirmative answers.

\smallskip Our main result in this paper, Theorem \ref{ikar}, immediately
and obviously implies that the answer to Question C is yes in one more
particular case, namely where $(A_{0},A_{1})$ is an arbitrary Banach couple,
and $(B_{0},B_{1})$ is the special couple of sequence spaces $(FL^{\infty
},FL_{1}^{\infty })$ of Fourier coefficients which was introduced and
studied in \cite{Ja}.

To the best of our knowledge, this particular case of Calder\'{o}n's problem
cannot be resolved via any other methods or cases treated in previous
partial solutions of the problem.

Our main motivation for considering this particular case is that a
corresponding affirmative (or negative!) result for the apparently related
couple $(B_{0},B_{1})=(\ell ^{\infty }(FL^{\infty }),\ell ^{\infty
}(FL_{1}^{\infty }))$ would resolve Calder\'{o}n's problem in complete
generality. (See \cite{CP}, \cite{Cw2} and \cite{CKM}.) It would apparently
also be sufficient to consider the couple of spaces of 
distributions $(B_{0},B_{1})=$ $(FL^{\infty }(\Bbb{R)},FL_{1}^{\infty }(\Bbb{%
R)})$. As shown in \cite{CKM}, it would also be sufficient to consider the
couple $(B_{0},B_{1})=(\ell ^{\infty }(FC),\ell ^{\infty }(FC_{1}))$.

\smallskip We believe that there ought to be a simpler proof of Theorem \ref
{ikar} than the one given here. One of the main underlying ideas is the use
of Lusin's theorem, which leads us to the intermediate result Theorem \ref
{diag}. Intuitively it seems rather clear that Theorem \ref{diag} implies
Theorem \ref{ikar}. But, unfortunately, we have not found any way so far to
bypass the lengthy explanation which we need for the technical details of
this implication.

\smallskip

We thank Fedor Nazarov for many interesting discussions and enlightening
comments about Question C in general.

\section{Some preliminaries}

\smallskip

As usual we let $\Bbb{T}$ denote the torus $\{z\in \Bbb{C}:\left| z\right|
=1\}$. Define the $\sigma $-algebra $\mathcal{T}$ of measurable subsets of $%
\Bbb{T}$ to consist of all images of Lebesgue measurable subsets of $[0,2\pi
)$ under the bijection $\phi (t)=e^{it}$ and define $\mu $ to be the measure
on these sets defined by taking $\mu (E)$ to be the Lebesgue measure of $%
\phi ^{-1}(E)$. In other words, $\mu $ is simply arc length measure on $\Bbb{%
T}$.

The \textit{distance} $d(z_{1},z_{2})$ between two points $z_{1}$ and $z_{2}$
of $\Bbb{T}$ will be taken to be the length of the shortest arc in $\Bbb{T}$
joining them, i.e., $d(z_{1},z_{2})=d(1,z_{2}/z_{1})=|t|$ where $t$ is the
unique number in $(-\pi ,\pi ]$ which satisfies $z_{2}/z_{1}=e^{it}$.

For each $z\in \Bbb{T}$ and each $r\in (0,\pi ]$, we let $\Gamma (z,r)$ be
the closed arc in $\Bbb{T}$ of length $2r$ centred at $z$, i.e., $\Gamma
(z,r)=\left\{ ze^{it}:\left| t\right| \le r\right\} =\left\{ w\in \Bbb{T}%
:d(w,z)\le r\right\} $.

\smallskip Let $L^{1}(\Bbb{T)=}L^{1}(\Bbb{T},\mathcal{T},\mu \Bbb{)}$ and $%
L^{\infty }(\Bbb{T)=}L^{\infty }(\Bbb{T},\mathcal{T},\mu \Bbb{)}$ denote the
usual Lebesgue spaces of (equivalence classes of) measurable functions $f:%
\Bbb{T\rightarrow \Bbb{C}}$ normed by $\left\| f\right\| _{L^{1}}=\int_{\Bbb{%
T}}\left| f\right| d\mu =\int_{0}^{2\pi }\left| f(e^{it})\right| dt$ and $%
\left\| f\right\| _{L^{\infty }}=\mathrm{ess~sup}_{t\in [0,2\pi )}\left|
f(e^{it})\right| $ respectively. For each $r\in (0,\pi ]$ and $z\in \Bbb{T}$
the formula 
\begin{equation}
\int_{\Gamma (z,r)}fd\mu =\int_{-r}^{r}f(ze^{it})dt  \label{thronk}
\end{equation}
obviously holds whenever $f=\chi _{E}$ and $E\in \mathcal{T}$, and therefore
it also holds for all $f\in L^{1}(\Bbb{T)}$.

We conclude this section by recalling some standard facts related to the
Lebesgue differentiation theorem, and expressing them in ways that will be
convenient for our purposes here.

\begin{definition}
\label{omf}For each $\mu $-integrable function $f:\Bbb{T\rightarrow \Bbb{C}}$%
, let $\Omega _{f}$ denote the set of all points $z\in \Bbb{T}$ for which
the limit $\lim_{r\rightarrow 0}\frac{1}{2r}\int_{-r}^{r}f(ze^{it})dt$
exists. Let $f^{\circ }:\Bbb{T\rightarrow \Bbb{C}}$ denote the function
defined by $f^{\circ }(z)=\lim_{r\rightarrow 0}\frac{1}{2r}%
\int_{-r}^{r}f(ze^{it})dt\cdot \chi _{\Omega _{f}}(z)$ for all $z\in \Bbb{T}$%
. It will be convenient to refer to the points of $\Omega _{f}$ as the 
\textit{weak canonical Lebesgue points }of $f$.
\end{definition}

If two $\mu $-integrable functions $f$ and $g$ coincide almost everywhere on 
$\Bbb{T}$, then of course $\Omega _{f}=\Omega _{g}$ and so $f^{\circ
}(z)=g^{\circ }(z)$ for \textit{all} $z\in \Bbb{T}$. Thus our definitions of 
$\Omega _{f}$ and $f^{\circ }$ extend unambiguously to the case where $f$ is
an element of $L^{1}(\Bbb{T})$ rather than a single function. It will be
convenient to refer to the function{\LARGE \ }$f^{\circ }$ as the \textit{%
canonical representative} of the element $f$. It is an obvious consequence
of the Lebesgue differentiation theorem that $\Omega _{f}\in \mathcal{T}$
with $\mu (\Bbb{T}\backslash \Omega _{f})=0$ and that $f(z)=f^{\circ }(z)$
for almost every $z\in \Bbb{T}$. This in turn implies that $\mathrm{ess~sup}%
_{t\in [0,2\pi )}\left| f(e^{it})\right| =\mathrm{sup}_{t\in [0,2\pi
)}\left| f^{\circ }(e^{it})\right| $.

\begin{definition}
\label{wdp}For each set $E\in \mathcal{T}$, let $E^{(d)}$ denote the set of
all \textit{points of density} of $E$, i.e., the points $z$ in $\Bbb{T}$ for
which the limit $\lim_{r\rightarrow 0}\frac{1}{2r}\int_{-r}^{r}\chi
_{E}(ze^{it})dt$ exists and is equal to $1$, i.e.,\ 
\begin{equation*}
\lim_{r\rightarrow 0}\frac{\mu (E\cap \Gamma (z,r))}{2r}=1.
\end{equation*}
\end{definition}

\smallskip It is of course a classical result, deduced by applying the
Lebesgue differentiation theorem to the function $\chi _{E}$, that $%
E^{(d)}\in \mathcal{T}$ and $\mu (E^{(d)}\backslash E)=\mu (E\backslash
E^{(d)})=0$. Since $E=(E\cap E^{(d)})\cup E\backslash E^{(d)}$, this implies
in turn that $\mu (E\cap E^{(d)})=\mu (E)$ and so, of course, 
\begin{equation}
\mu (\Bbb{T}\backslash (E\cap E^{(d)}))=\mu (\Bbb{T}\backslash E).
\label{gurp}
\end{equation}

\section{The main ingredient}

Our goal in this section is to prove Theorem \ref{diag} which will
subsequently be the principal tool in the proof of our main result. But
first we need the following lemma:

\begin{lemma}
\label{qaa}Let $K$ be a relatively compact subset of $L^{\infty }(\Bbb{T})$.
Then, for each $\epsilon >0$, there exists a number $\delta =\delta
(\epsilon )$ and a set $U=U_{\epsilon }\subset \Bbb{T}$ such that 
\begin{equation}
\mu (\Bbb{T}\backslash U)<\epsilon   \label{urp}
\end{equation}
and, for each $f\in K$, 
\begin{equation*}
\mathrm{ess~sup}\left\{ \left| f(z_{1})-f(z_{2})\right| :z_{1},z_{2}\in
U,d(z_{1},z_{2})<\delta \right\} \le \epsilon .
\end{equation*}
More precisely, 
\begin{equation}
\sup \left\{ \left| f^{\circ }(z_{1})-f^{\circ }(z_{2})\right|
:z_{1},z_{2}\in U\cap \Omega _{f},d(z_{1},z_{2})<\delta \right\} \le
\epsilon .  \label{vrim}
\end{equation}
\end{lemma}

\smallskip

\begin{proof} 
Since $K$ is relatively compact, there exists some integer $%
N $ and a collection of $N$ elements $\left\{ f_{1},f_{2},...,f_{N}\right\}
\subset K$, such that $K\subset \bigcup_{m=1}^{N}B(f_{j},\epsilon /3)$.
(Here of course $B(f,r)$ denotes the set $\{g\in L^{\infty }(\Bbb{T)}%
:\left\| f-g\right\| _{L^{\infty }(\Bbb{T)}}<r\}$.)

The $\sigma $-algebra $\mathcal{T}$ contains all the Borel sets in $\Bbb{T}$%
, and the measure $\mu $ is finite, regular and complete on $\mathcal{T}$,
i.e., it satisfies conditions (b), (c) and (d) of Theorem 2.14 of \cite
{Rudin} pp.\ 40--41. Thus we can apply Lusin's Theorem (\cite{Rudin} Theorem
2.24 p.\ 55) to obtain, for each $j\in \{1,2,...,N\}$, a continuous function 
$g_{j}:\Bbb{T}\rightarrow \Bbb{C}$ such that the set $G_{j}=\{z\in \Bbb{T}%
:f_{j}^{\circ }(z)\ne g_{j}(z)\}$ has measure $\mu (G_{j})<\epsilon /N$.
Since each $g_{j}$ is uniformly continuous, we can choose $\delta =\delta
(\epsilon )$ to be a number with the property that 
\begin{equation*}
z_{1},z_{2}\in \Bbb{T},d(z_{1},z_{2})<\delta \Longrightarrow \max \left\{
\left| g_{j}(z_{1})-g_{j}(z_{2})\right| :j=1,2,...,N\right\} <\epsilon /3.
\end{equation*}

Let $E$ be the set $E=\Bbb{T}\backslash \bigcup_{j=1}^{N}G_{j}$\smallskip .
Then we will choose the set $U$ to be $U=E\cap E^{(d)}$. Clearly $\Bbb{T}%
\backslash E=\bigcup_{j=1}^{N}G_{j}$ has measure $\mu \left(
\bigcup_{j=1}^{N}G_{j}\right) $ not exceeding $\sum_{j=1}^{N}\mu
(G_{j})<\epsilon $, and so (\ref{urp}) follows immediately from (\ref{gurp}).

\smallskip For each point $z\in \Bbb{T}$ \smallskip and each $r\in (0,\pi ]$
and each $j\in \{1,2,...,N\}$, we have 
\begin{eqnarray*}
&&\int_{\Gamma (z,r)}\left| f_{j}^{\circ }(w)-g_{j}(w)\right| d\mu (w) \\
&=&\int_{\Gamma (z,r)\cap E}\left| f_{j}^{\circ }(w)-g_{j}(w)\right| d\mu
(w)+\int_{\Gamma (z,r)\backslash E}\left| f_{j}^{\circ }(w)-g_{j}(w)\right|
d\mu (w) \\
&=&0+\int_{\Gamma (z,r)\backslash E}\left| f_{j}^{\circ }(w)-g_{j}(w)\right|
d\mu (w) \\
&\le &\mu (\Gamma (z,r)\backslash E)\left( \left\| f_{j}\right\| _{L^{\infty
}}+\left\| g_{j}\right\| _{L^{\infty }}\right)
\end{eqnarray*}

If $z\in U$, and is therefore a point of density of $E$, then the preceding
estimates show that 
\begin{multline*}
\limsup_{r\rightarrow 0}\frac{1}{\mu (\Gamma (z,r))}\int_{\Gamma (z,r)}\left|
f_{j}^{\circ }(w)-g_{j}(w)\right| d\mu (w)
\\
\le \left( \left\| f_{j}\right\|
_{L^{\infty }}+\left\| g_{j}\right\| _{L^{\infty }}\right)
\lim_{r\rightarrow 0}\frac{\mu (\Gamma (z,r))-\mu ((\Gamma (z,r)\cap E))}{%
\mu (\Gamma (z,r))}=0.
\end{multline*}

In other words (cf.\ (\ref{thronk}), we have 
\begin{equation}
\lim_{r\rightarrow 0}\frac{1}{2r}\int_{-r}^{r}\left| f_{j}^{\circ
}(ze^{it})-g_{j}(ze^{it})\right| dt=0
\quad
\text{for all }z\in U\text{ and all }%
j\in \{1,2,...,N\}.  \label{zblurk}
\end{equation}

\smallskip

Now suppose that $f$ is some arbitrary element of $K$. Then there exists
some integer $j$ depending on $f$ such that 
\begin{equation}
f\in B(f_{j},\epsilon /3).  \label{klonk}
\end{equation}

Suppose that $z_{1}$ and $z_{2}$ are arbitrary points in $U\cap \Omega _{f}$
satisfying $d(z_{1},z_{2})<\delta $. For each $w\in \Bbb{T}$ we have 
\begin{multline}
\left| f^{\circ }(wz_{1})-f^{\circ }(wz_{2})\right|  
\le \left| f^{\circ }(wz_{1})-f_{j}^{\circ }(wz_{1})\right| +\left|
f_{j}^{\circ }(wz_{1})-g_{j}(wz_{1})\right| \\+\left|
g_{j}(wz_{1})-g_{j}(wz_{2})\right|  \label{zjonj} 
+\left| g_{j}(wz_{2})-f_{j}^{\circ }(wz_{2})\right| +\left| f_{j}^{\circ
}(wz_{2})-f^{\circ }(wz_{2})\right| .  
\end{multline}
Since $d(wz_{1},wz_{2})<\delta $, we have $\left|
g_{j}(wz_{1})-g_{j}(wz_{2})\right| <\epsilon /3$ for each $w\in \Bbb{T}$. In
view of (\ref{klonk}) we also have $\left| f^{\circ }(wz_{1})-f_{j}^{\circ
}(wz_{1})\right| <\epsilon /3$ and 
$\left| f^\circ(wz_{2})-f^\circ_{j}(wz_{2})\right|
<\epsilon /3$ for  every $w\in \Bbb{T}$. So, for each $r\in (0,\pi ]$,
we can set $w=e^{it}$ and integrate the inequality (\ref{zjonj}) with
respect to $t$ on the interval $[-r,r]$ to obtain that 
\begingroup
\multlinegap=0pt
\begin{multline*}
\frac{1}{2r}\int_{-r}^{r}\left| f^{\circ }(e^{it}z_{1})-f^{\circ
}(e^{it}z_{2})\right| dt \\
\le \epsilon /3+\frac{1}{2r}\int_{-r}^{r}\left| f_{j}^{\circ
}(e^{it}z_{1})-g_{j}(e^{it}z_{1})\right| dt+\epsilon /3+\frac{1}{2r}%
\int_{-r}^{r}\left| g_{j}(e^{it}z_{2})-f_{j}^{\circ }(e^{it}z_{2})\right|
dt+\epsilon /3.
\end{multline*}
\endgroup
In view of (\ref{zblurk}) both of the integrals in the preceding line tend
to $0$ as $r$ tends to $0$, and so we have shown that 
\begin{equation}
\limsup_{r\rightarrow 0}\frac{1}{2r}\int_{-r}^{r}\left|
f(e^{it}z_{1})-f(e^{it}z_{2})\right| dt\le \epsilon .  \label{obh}
\end{equation}
Finally, since $z_{1}$ and $z_{2}$ are both in $\Omega _{f}$, we see that 
\begin{eqnarray*}
\left| f^{\circ }(z_{1})-f^{\circ }(z_{2})\right| &=&\left|
\lim_{r\rightarrow 0}\frac{1}{2r}\int_{-r}^{r}f(e^{it}z_{1})dt-\lim_{r%
\rightarrow 0}\frac{1}{2r}\int_{-r}^{r}f(e^{it}z_{2})dt\right| \\
&=&\lim_{r\rightarrow 0}\left| \frac{1}{2r}
\int_{-r}^{r}\bigl(f(e^{it}z_{1})-f(e^{it}z_{2})\bigr)dt\right| \\
&\le &\limsup_{r\rightarrow 0}\frac{1}{2r}\int_{-r}^{r}\left|
f(e^{it}z_{1})-f(e^{it}z_{2})\right| dt.
\end{eqnarray*}

This, combined with (\ref{obh}) establishes (\ref{vrim}) and completes the
proof of Lemma \ref{qaa}. 
\end{proof}

\smallskip

In order to be able to conveniently apply the next theorem later in the
proof of our main result, we will have to distinguish, more pedantically
than is usually necessary, between elements of $L^{\infty }(\Bbb{T})$ and
the functions that represent them. So here we will use bold lower case
letters for elements (equivalent classes of complex valued functions) in $%
L^{\infty }(\Bbb{T})$ and usual italic lower case letters for the complex
valued functions in those equivalence classes.

\begin{theorem}
\smallskip \label{diag}Let $K$ be a subset of $L^{\infty }(\Bbb{T)}$. Let $%
\mathcal{E}_{K}$ be the set of all functions $u:\Bbb{T}\times \Bbb{T}%
\rightarrow \Bbb{C}$ which are finite sums of the form 
\begin{equation}
u(\zeta ,w)=\sum_{n=-M}^{M}\zeta ^{n}v_{n}(w),  \label{sotf}
\end{equation}
where each function $v_{n}:\Bbb{T}\rightarrow \Bbb{C}$ is essentially
bounded and where the elements $\mathbf{v}_{n}\in L^{\infty }(\Bbb{T})$,
which are represented by $v_{n}$ respectively for each integer $n\in [-M,M]$%
, satisfy 
\begin{equation}
\sum_{n=-M}^{M}\zeta ^{n}\mathbf{v}_{n}\in K\quad\text{for each }\zeta
\in \Bbb{T}.  \label{gbp}
\end{equation}

For each $\delta >0$, define 
\begin{equation}
\rho _{K}(\delta ):=\sup \left\{ \int_{\Bbb{T}}\left|
u(w,z_{1}w)-u(w,z_{2}w)\right| d\mu (w):u\in \mathcal{E}_{K},z_{1},z_{2}\in 
\Bbb{T},d(z_{1},z_{2})\le \delta \right\} .  \label{droq}
\end{equation}
If $K$ is a relatively compact subset of $L^{\infty }(\Bbb{T})$ then 
\begin{equation}
\lim_{\delta \rightarrow 0}\rho _{K}(\delta )=0.  \label{ucay}
\end{equation}
\end{theorem}

\begin{proof}
Let us begin by checking the almost obvious fact that the
integral \linebreakx
$\int_{\Bbb{T}}\left| u(w,z_{1}w)-u(w,z_{2}w)\right| d\mu (w)$ does
not depend on our particular choices of the representatives $v_{n}$ of $%
\mathbf{v}_{n}$ in the formula (\ref{sotf}).

\smallskip More explicitly, we claim that if $\phi (\zeta
,w)=\sum_{n=-M}^{M}\zeta ^{n}f_{n}(w)$ where $f_{n}(w)=v_{n}(w)$ for a.e.\ $%
w\in \Bbb{T}$ and for each $n$, then 
\begin{equation}
\int_{\Bbb{T}}\left| u(w,z_{1}w)-u(w,z_{2}w)\right| d\mu (w)=\int_{\Bbb{T}%
}\left| \phi (w,z_{1}w)-\phi (w,z_{2}w)\right| d\mu (w).  \label{ppz}
\end{equation}

To establish (\ref{ppz}) we first observe that, for each fixed $z\in
\Bbb{T}$, 
\begin{eqnarray}
\int_{\Bbb{T}}\left| \phi (w,zw)-u(w,zw)\right| d\mu (w) &=&\int_{\Bbb{T}%
}\left| \sum_{n=-M}^{M}w^{n}(f_{n}(zw)-v_{n}(zw))\right| d\mu (w)  \notag \\
&\le &\sum_{n=-M}^{M}\int_{\Bbb{T}}\left| w^{n}(f_{n}(zw)-v_{n}(zw))\right|
d\mu (w)  \notag \\
&=&\sum_{n=-M}^{M}\int_{\Bbb{T}}\left| f_{n}(zw)-v_{n}(zw)\right| d\mu (w) 
\notag \\
&=&\sum_{n=-M}^{M}\int_{\Bbb{T}}\left| f_{n}(w)-v_{n}(w)\right| d\mu (w)=0.
\label{iofg}
\end{eqnarray}

\smallskip Next we note that, for each fixed $z_{1}$ and $z_{2}$, 
\begin{multline*}
\left| \phi (w,z_{1}w)-\phi (w,z_{2}w)\right| \le \left| \phi
(w,z_{1}w)-u(w,z_{1}w)\right| +
\\
\left| u(w,z_{1}w)-u(w,z_{2}w)\right| +\left|
u(w,z_{2}w)-\phi (w,z_{2}w)\right| .
\end{multline*}
Integrating, and using (\ref{iofg}), we obtain that 
\begin{equation}
\int_{\Bbb{T}}\left| \phi (w,z_{1}w)-\phi (w,z_{2}w)\right| d\mu (w)\le
\int_{\Bbb{T}}\left| u(w,z_{1}w)-u(w,z_{2}w)\right| d\mu (w).  \label{gomh}
\end{equation}
The reverse inequality to (\ref{gomh}) is obtained exactly analogously and
so indeed we obtain (\ref{ppz}).

\smallskip We will use the usual notation $\eta G=\left\{ \eta w:w\in
G\right\} $ for each $\eta \in \Bbb{C}$ and for each subset $G$ of some
complex vector space $W$ (which will often also be simply $\Bbb{\Bbb{C}}$).

\smallskip

Obviously the sets $\mathcal{E}_{K}$ and the quantities $\rho _{K}(\delta )$
satisfy $\mathcal{E}_{\eta K}=\eta \mathcal{E}_{K}$ and so $\rho _{\eta
K}(\delta )=\eta \rho _{K}(\delta )$ for each $\eta >0$. So we can suppose,
without loss of generality, that the given relatively compact set $K$ is
contained in the unit ball of $L^{\infty }(\Bbb{T})$.

\smallskip

Choose a positive number $\epsilon $ and let $U=U_{\epsilon }$ and $\delta
=\delta (\epsilon )$ be a set and a number with the properties listed in the
statement of Lemma \ref{qaa}, corresponding to our current choices of $%
\epsilon $ and of the relatively compact set $K$. We will show that, for
this choice of $\delta $, 
\begin{equation}
\sup \left\{ \int_{\Bbb{T}}\left| u(w,z_{1}w)-u(w,z_{2}w)\right| d\mu
(w):u\in \mathcal{E}_{K},z_{1},z_{2}\in \Bbb{T},d(z_{1},z_{2})\le \delta
\right\} \le (2\pi +4)\epsilon .  \label{ecay}
\end{equation}
Since $\epsilon $ can be chosen arbitrarily small, this will of course
suffice to establish (\ref{ucay}).

Let $u(\zeta ,w)=\sum_{n=-M}^{M}\zeta ^{n}v_{n}(w)$ be an arbitrary element
of $\mathcal{E}_{K}$. In view of (\ref{ppz}) we may assume without loss of
generality that it is in fact given by the formula 
\begin{equation}
u(\zeta ,w)=\sum_{n=-M}^{M}\zeta ^{n}v_{n}^{\circ }(w)  \label{tvh}
\end{equation}
where, for each $n$, $v_{n}^{\circ }$ is the canonical representative of $%
v_{n}$, and of $\mathbf{v}_{n}$ the element of $L^{\infty }(\Bbb{T})$
represented by $v_{n}$).

Let $\Omega =\bigcap_{n=-M}^{M}\Omega _{v_{n}}$. Clearly 
\begin{equation}
\mu (\Bbb{T}\backslash \Omega )=\mu \left( \bigcup_{n=-M}^{M}\Bbb{T}%
\backslash \Omega _{v_{n}}\right) =0,  \label{bbw}
\end{equation}
and, for each fixed $\zeta \in \Bbb{T}$ , every point in $\Omega $ is a weak
canonical Lebesgue point for the function $w\mapsto u(\zeta ,w)$. (Of course
the set of all weak canonical Lebesgue points for this function may be
strictly larger than $\Omega $.)

\smallskip In view of (\ref{gbp}) and the fact that $K$ is a subset of the
unit ball of $L^{\infty }(\Bbb{T})$, it follows that for each fixed $\zeta
\in \Bbb{T}$, we have $\left| u(\zeta ,w)\right| \le 1$ for almost every $%
w\in \Bbb{T}$. 
The exceptional null set may depend on $\zeta$, but the next simple lemma
shows that the exceptional null sets for different $\zeta$ are
contained in a common null set. 

\begin{lemma}\label{Lpddl}
Suppose that $u(\zeta ,w)=\sum_{n=-M}^{M}\zeta ^{n}v_{n}(w)$ for
some fixed integer $M$ and all $\zeta ,w\in \Bbb{T}$, 
where $v_n$ are some measurable functions on $\Bbb T$,
and suppose that, 
for each fixed $\zeta \in \Bbb{T}$, $\left| u(\zeta ,w)\right| \le 1$
for almost all $w\in \Bbb{T}$. 
Then, for almost all $w\in \Bbb{T}$, $\left| u(\zeta ,w)\right| \le 1$
holds for all $\zeta \in \Bbb{T}$. 
\end{lemma}

\begin{proof}
Let, for $\zeta\in\Bbb T$,  $N_\zeta$ be the null set
\set{w:\left| u(\zeta ,w)\right| > 1}. Let $Q$ be a countable dense
subset of $\Bbb T$ (for example, the set of points with rational argument),
and let $N=\bigcup_{\zeta\in Q} N_\zeta$. The $\mu(N)=0$, and if
$w\notin N$, then $\left| u(\zeta ,w)\right| \le 1$ for all 
$\zeta\in Q$. 
However, for fixed $w$, $u(\zeta,w)$ is a continuous function of
$\zeta$, so the inequality holds for all $\zeta\in\Bbb T$.
\end{proof}

In particular, for every $z\in\Bbb T$, 
$|u(w,zw)|\le1$ for a.e.\ $w$, and thus
\begin{equation}
\int_{F}\left| u(w,zw)\right| d\mu (w)\le \mu (F)
\quad\text{for every set }F\in \mathcal{T} \text{ and }z\in\Bbb T.  
\label{ddl}
\end{equation}

Obviously $\mu (\eta G)=\mu (G)$ for every measurable set $G\subset \Bbb{T}$
and every $\eta \in \Bbb{T}$. Let $z_{1}$ and $z_{2}$ be arbitrary points in 
$\Bbb{T}$ satisfying 
\begin{equation}
d(z_{1},z_{2})<\delta  \label{ammm}
\end{equation}
and consider the set $F=$ $z_{1}^{-1}(U\cap \Omega )\cap z_{2}^{-1}(U\cap
\Omega )$. Whenever $w\in F$, both the points $z_{1}w$ and $z_{2}w$ are in $%
U\cap \Omega $. In view of (\ref{ammm}) we also have $d(z_{1}w,z_{2}w)<%
\delta $. Furthermore, (by (\ref{gbp}) and (\ref{tvh})), for any fixed value
of $w$ in $F$, the function $f(z):=u(w,z)$ is a representative of an element
of $K$, the points $z_{1}w$ and $z_{2}w$ are both in $\Omega _{f}$, 
and $f^\circ(z_1w)=f(z_1w)$, $f^\circ(z_2w)=f(z_2w)$.
Consequently, applying (\ref{vrim}), we have 
\begin{equation}
\left| u(w,z_{1}w)-u(w,z_{2}w)\right| =\left| f(z_{1}w)-f(z_{2}w)\right| \le
\epsilon \quad\text{for all }w\in F.  \label{mpcm}
\end{equation}

Now we need to estimate the measure of the complement of the set $F$. We
write $F=F_{1}\cap F_{2}$, where $F_{j}=z_{j}^{-1}(U\cap \Omega )$ for $%
j=1,2 $ First, since $U=(U\cap \Omega )\cup (U\backslash \Omega )$, we
deduce from (\ref{bbw}) that $\mu (U\cap \Omega )=\mu (U)$. Therefore, $\mu
(F_{1})=\mu (F_{2})=\mu (U)$ and $\mu (\Bbb{T}\backslash F_{1})=\mu (\Bbb{T}%
\backslash F_{2})=\mu (\Bbb{T}\backslash U)$. We also recall that the
properties listed in the statement of Lemma \ref{qaa}, which $U$ must
satisfy, include the condition $\mu (\Bbb{T}\backslash U)<\epsilon $. From
these remarks, and the fact that $\Bbb{T}\backslash (F_{1}\cap F_{1})=\left( 
\Bbb{T}\backslash F_{1}\right) \cup (\Bbb{T}\backslash F_{2})$, we can now
see that 
\begin{equation}
\mu \left( \Bbb{T}\backslash F\right) <2\epsilon .  \label{pst}
\end{equation}

\smallskip Finally, by (\ref{mpcm}), (\ref{ddl}) and (\ref{pst}), we obtain
that
\begin{eqnarray*}
&&\int_{\Bbb{T}}\left| u(w,z_{1}w)-u(w,z_{2}w)\right| d\mu (w) \\
&=&\int_{F}\left| u(w,z_{1}w)-u(w,z_{2}w)\right| d\mu (w)+\int_{\Bbb{T}%
\backslash F}\left| u(w,z_{1}w)-u(w,z_{2}w)\right| d\mu (w) \\
&\le &\int_{F}\epsilon d\mu (w)+\int_{\Bbb{T}\backslash F}\left|
u(w,z_{1}w)\right| d\mu (w)+\int_{\Bbb{T}\backslash F}\left|
u(w,z_{2}w)\right| d\mu (w) \\
&\le &\mu (F)\epsilon +2\mu (\Bbb{T}\backslash F)\le (2\pi +4)\epsilon .
\end{eqnarray*}

Since the element $u$ of $\mathcal{E}_{K}$ and the points $z_{1}$ and $z_{2}$
in $\Bbb{T}$ satisfying $d(z_{1},z_{2})<\delta $ were chosen arbitrarily,
this establishes (\ref{ecay}) and so completes the proof of Theorem~\ref
{diag}. 
\end{proof}

\smallskip

\section{Some more preliminaries}

Next we recall some standard things about Alberto Calder\'{o}n's complex
interpolation spaces. We choose versions of the notation and definitions to
suit our particular purposes. Some of the facts that we prove are implicit
or even sometimes explicit in Calder\'{o}n's fundamental paper \cite{Ca}
about these spaces, or in other papers. But we prefer to give rather full
explanations here. We will assume some familiarity with the basic notions
and terminology of interpolation space theory as presented, for example, in
the early chapters of \cite{bl} or of \cite{bk}.

\smallskip Let $\Bbb{S}$ be the ``unit strip'', i.e., $\Bbb{S}=\left\{ z\in 
\Bbb{C}:0\le \mathrm{Re}\,z\le 1\right\} $ and let $\Bbb{A}$ denote the
``unit annulus'', $\Bbb{A}=\left\{ z\in \Bbb{C}:1\le \left| z\right| \le
e\right\} $. Calder\'{o}n constructed his complex interpolation spaces in 
\cite{Ca} via certain analytic vector valued functions defined on $\Bbb{S}$.
Most of the time here we prefer to use an alternative construction (cf.\ \cite
{Cw1}) where $\Bbb{S}$ is replaced by $\Bbb{A}$. 

Let $\vec{B}=(B_{0},B_{1})$ be an arbitrarily Banach couple of complex
Banach spaces. Let $\mathcal{F}_{\Bbb{A}}(\vec{B})=\mathcal{F}_{\Bbb{A}%
}(B_{0},B_{1})$ be the space of all continuous functions $f:\Bbb{A}%
\rightarrow B_{0}+B_{1}$ such that
\begin{equation}
f\text{ is analytic in the interior }\Bbb{A}^{\circ }\text{ of }\Bbb{A}
\label{ffto}
\end{equation}
and 
\begin{equation}
\text{for }j=0,1\text{, the restriction of }f\text{ to the circle }e^{j}\Bbb{%
T}\text{ is a continuous map of }e^{j}\Bbb{T}\text{ into }B_{j}\text{ .}
\label{fftt}
\end{equation}

We norm $\mathcal{F}_{\Bbb{A}}(\vec{B})$ by $\left\| f\right\| _{\mathcal{F}%
_{\Bbb{A}}(\vec{B})}=\sup \left\{ \left\| f(e^{j+it})\right\| _{B_{j}}:t\in
[0,2\pi ),j=0,1\right\} $ .

\smallskip For each $\theta \in (0,1)$, let $[\vec{B}]_{\theta ,\Bbb{A}}$
denote the space of all elements in $B_{0}+B_{1}$ of the form $b=f(e^{\theta
})$, for some $f\in \mathcal{F}_{\Bbb{A}}(\vec{B})$. It is normed by 
$$\left\| b\right\| _{[\vec{B}]_{\theta ,\Bbb{A}}}:=\inf \left\{ \left\|
f\right\| _{\mathcal{F}_{\Bbb{A}}(\vec{B})}:f\in \mathcal{F}_{\Bbb{A}}(\vec{B%
}),f(e^{\theta })=b\right\} .$$  
As shown in \cite{Cw1}, the space $%
[\vec{B}]_{\theta ,\Bbb{A}}$, coincides, to within equivalence of norms,
with Calder\'{o}n's complex interpolation space $[\vec{B}]_{\theta }$.

\smallskip Each function $f\in \mathcal{F}_{\Bbb{A}}(\vec{B})$, can be
expressed as a Laurent series $f(z)=\sum_{n=-\infty }^{\infty }z^{n}\widehat{%
f}(n)$, which converges absolutely and uniformly in $B_{0}+B_{1}$ norm on
every compact subset of $\Bbb{A}^{\circ }$. The coefficients $\widehat{f}(n)$
are given by the formul\ae\ 
\begin{equation}
\widehat{f}(n)=\frac{1}{2\pi i}\oint_{\beta \Bbb{T}}\frac{1}{z^{n+1}}f(z)dz,
\label{fflc}
\end{equation}
where here, and in what follows, (cf.\ \cite{Ka} p.\ 10 and pp.\ 257--8) we
are using Riemann integrals of continuous Banach space valued functions.
Here $\beta $ is a number in the interval $(1,e)$ and the value of $\widehat{%
f}(n)$ is of course independent of its particular value in $(1,e)$. We can
of course rewrite (\ref{fflc}) as 
\begin{equation}
\widehat{f}(n)=\frac{1}{2\pi }\int_{0}^{2\pi }e^{-n(\alpha +it)}f(e^{\alpha
+it})dt\text{ where }\alpha =\ln \beta \in (0,1).  \label{2fflc}
\end{equation}

For each fixed $n$ the function $\alpha \mapsto \frac{1}{2\pi }%
\int_{0}^{2\pi }e^{-n(\alpha +it)}f(e^{\alpha +it})dt$ is also defined for $%
\alpha =0$ and $\alpha =1$, and is a continuous map of the closed interval $%
[0,1]$ into $B_{0}+B_{1}$. Consequently, the formula (\ref{2fflc}) also
holds for $\alpha =0$ and $\alpha =1$. For these two values, the integral $%
\frac{1}{2\pi }\int_{0}^{2\pi }e^{-n(\alpha +it)}f(e^{\alpha +it})dt$
defines an element of $B_{0}$ or of $B_{1}$ respectively. Thus we see that $%
\widehat{f}(n)\in B_{0}\cap B_{1}$ for each $n\in \Bbb{Z}$.

For each $f\in \mathcal{F}_{\Bbb{A}}(\vec{B})$ and each $N\in \Bbb{N}$, we
define the $N$th Fej\'{e}r mean of $f$ to be the function $\sigma _{N}(f):%
\Bbb{A}\rightarrow B_{0}+B_{1}$ given by the formula 
\begin{equation*}
\sigma _{N}(f)(e^{\alpha +it})=\frac{1}{2\pi }\int_{0}^{2\pi }\kappa
_{N}(t-s)f(e^{\alpha +is})ds\quad\text{for each }\alpha \in [0,1]\text{ and }%
t\in [0,2\pi ),
\end{equation*}
where $\kappa _{N}$ is the usual Fej\'{e}r kernel, $\kappa _{N}(t)=\frac{1}{%
N+1}\left( \frac{\sin \frac{N+1}{2}t}{\sin \frac{t}{2}}\right) ^{2}$. Since
(cf.\ \cite{Ka} p.\ 12) $\kappa _{N}(t)=\sum_{n=-N}^{N}\left( 1-\frac{\left|
n\right| }{N+1}\right) e^{int}$, it follows immediately that $\ \sigma
_{N}(f)(e^{\alpha +it})=\sum_{n=-N}^{N}\left( 1-\frac{\left| n\right| }{N+1}%
\right) e^{n(\alpha +it)}\widehat{f}(n)$, i.e., $\sigma
_{N}(f)(z)=\sum_{n=-N}^{N}\left( 1-\frac{\left| n\right| }{N+1}\right) z^{n}%
\widehat{f}(n)$ for all $z\in \Bbb{A}$. It follows obviously from (\ref{fftt}%
) that $\left. f\right| _{e^{j}\Bbb{T}}:e^{j}\Bbb{T\rightarrow }B_{j}$ is in
fact uniformly continuous for $j=0,1$. This combined with the fact that $%
\left\{ \kappa _{N}\right\} _{N\in \Bbb{N}}$ is a summability kernel (\cite
{Ka} p.\ 9) readily implies, using standard estimates like those in the
proof of \cite{Ka} Lemma 2.2 p.\ 10, that $\lim_{N\rightarrow \infty }\sup
\left\{ \left\| \sigma _{N}(f)(e^{j+it})-f(e^{j+it})\right\| _{B_{j}}:t\in
[0,2\pi )\right\} =0$ for $j=0,1$. Thus we have that $\lim_{N\rightarrow
\infty }\left\| \sigma _{N}(f)-f\right\| _{\mathcal{F}_{\Bbb{A}}(\vec{B})}=0$
and consequently also 
\begin{equation}
\lim_{N\rightarrow \infty }\left\| \sigma _{N}(f)(e^{\theta })-f(e^{\theta
})\right\| _{[\vec{B}]_{\theta ,\Bbb{A}}}=0\quad\text{for each }\theta
\in (0,1). 
\label{pgdensef}
\end{equation}

Let $\mathcal{G}_{\Bbb{A}}(\vec{B})$ denote the subspace of $\mathcal{F}_{%
\Bbb{A}}(\vec{B})$ consisting of all functions $g:\Bbb{A}\rightarrow
B_{0}\cap B_{1}$ which are finite sums of the form $g(z)=%
\sum_{n=-N}^{N}z^{n}b_{n}$, for elements $b_{n}\in B_{0}\cap B_{1}$. The
remarks in the preceding paragraph show that 
$\sigma_N(f)\in\mathcal{G}_{\Bbb{A}}(\vec{B})$
for every $f\in\mathcal{F}_{\Bbb{A}}(\vec{B})$, and that hence
$\mathcal{G}_{\Bbb{A}}(\vec{B})$
is dense in $\mathcal{F}_{\Bbb{A}}(\vec{B})$.

\begin{definition}
\label{fkk}For any given subsets $K_{0}$ of $B_{0}$ and $K_{1}\subset B_{1}$%
, we let $\mathcal{F}_{\Bbb{A}}(K_{0},K_{1})$ denote the subset of functions 
$f\in \mathcal{F}_{\Bbb{A}}(\vec{B})$ with the additional property that $%
f(e^{j+it})\in K_{j}$ for all $t\in [0,2\pi )$ and $j=0,1$.
\end{definition}

\begin{fact}
\label{konv}Let $K_{j}$ be an arbitrary subset of $B_{j}$ for $j=0,1$. Let $%
\widetilde{K_{j}}$ denote the closure in $B_{j}$ of the convex hull of $K_{j}
$ and let $K_{j}^{*}$ denote the closure in $B_{j}$ of the absolutely convex
hull of $K_{j}$. Let $f$ be an arbitrary element of $\mathcal{F}_{\Bbb{A}%
}(K_{0},K_{1})$. Then
\begin{romenumerate}
\item $\widehat{f}(n)\in K_{0}^{*}\cap K_{1}^{*}$ for each $n\in \Bbb{Z}$, and

\item 
the function $\sigma _{N}(f)$ is in $\mathcal{F}_{\Bbb{A}}(\widetilde{%
K_{0}},\widetilde{K_{1}})$ for each $N\in \Bbb{N}$.
\end{romenumerate}
\end{fact}

\begin{proof}[Proof of Fact \ref{konv}]
Perhaps the quickest and easiest way to
see (ii) is to consider the function $f_{M}:\Bbb{A\rightarrow }B_{0}$
defined by,
for each $t\in [0,2\pi )$ and each $\alpha \in [0,1]$. 
\begin{equation*}
f_{M}(e^{\alpha +it})=\frac{1}{2\pi }\sum_{m=1}^{M}\int_{2\pi (m-1)/M}^{2\pi
m/M}\kappa _{N}(t-s)f(e^{\alpha +i2\pi m/M})\,ds.
%\quad\text{for each }t\in [0,2\pi ) \text{ and each }\alpha \in [0,1].
\end{equation*}
Since $\kappa _{N}\ge 0$ and $\frac{1}{2\pi }\sum_{m=1}^{M}\int_{2\pi
(m-1)/M}^{2\pi m/M}\kappa _{N}(t-s)ds=1$, it is clear that, for $j=0,1$, the
element $f_{M}(e^{j+it})$ is a convex combination of elements of $K_{j}$.
Furthermore, the uniform continuity of $\left. f\right| _{e^{j}\Bbb{T}}:e^{j}%
\Bbb{T\rightarrow }B_{j}$ ensures that $\lim_{M\rightarrow \infty }\left\|
f_{M}(e^{j+it})-\sigma _{N}(f)(e^{j+it})\right\| _{B_{0}}=0$. So indeed $%
\sigma _{N}(f)(e^{j+it})\in \widetilde{K_{j}}$.

\smallskip Analogously, for (i), for each fixed $n\in \Bbb{Z}$ and each $%
M\in \Bbb{N}$, we consider the element 
\begin{equation*}
f_{M}(n,j)=\frac{1}{2\pi }\sum_{m=1}^{M}\int_{2\pi (m-1)/M}^{2\pi
m/M}e^{-n(j+i2\pi m/M)}f(e^{j+i2\pi m/M})dt\quad\text{for each }j\in \{0,1\}.
\end{equation*}
Since $f_{M}(n,j)=\frac{1}{M}\sum_{m=1}^{M}e^{-n(j+i2\pi m/M)}f(e^{j+i2\pi
m/M})$, it is clear that $f_{M}(n,j)\in K_{j}^{*}$ for $j=0,1$. The uniform
continuity of the function $g:\Bbb{T\rightarrow }B_{j}$ defined by $%
g(e^{it})=e^{-n(j+it)}f(e^{j+it})$ ensures that $\lim_{M\rightarrow \infty
}\left\| f_{M}(n,j)-\frac{1}{2\pi }\int_{0}^{2\pi
}e^{-n(j+it)}f(e^{j+it})dt\right\| _{B_{j}}=0$. Since $\frac{1}{2\pi }%
\int_{0}^{2\pi }e^{-n(j+it)}f(e^{j+it})dt=\widehat{f(}n)$ for $j=0$ and also
for $j=1$ we deduce that $\widehat{f}(n)\in K_{0}^{*}\cap K_{1}^{*}$. 
\end{proof}

\smallskip

\smallskip Let $FL^{1}$ be the space of complex sequences $\lambda =\left\{
\lambda _{n}\right\} _{n\in \Bbb{Z}}$ which are Fourier coefficients of
functions in $L^{1}(\Bbb{T)}$, i.e., such that, for some function $%
u=u_{\lambda }\in L^{1}(\Bbb{T)}$, 
\begin{equation*}
\lambda _{n}=\frac{1}{2\pi }\int_{0}^{2\pi }e^{-int}u_{\lambda }(e^{it})dt%
\quad\text{for all }n\in \Bbb{Z}.
\end{equation*}
We norm $FL^{1}$ by setting $\left\| \lambda \right\| _{FL^{\infty
}}=\left\| u_{\lambda }\right\| _{L^{1}(\Bbb{T)}}$. Analogously, we let $%
FL^{\infty }$ be the space of complex sequences $\lambda =\left\{ \lambda
_{n}\right\} _{n\in \Bbb{Z}}$ which are Fourier coefficients of functions in 
$L^{\infty }(\Bbb{T)}$, i.e., such that the above function $u_{\lambda }$ is
also in $L^{\infty }(\Bbb{T)}$. Here we use the norm $\left\| \lambda
\right\| _{FL^{\infty }}=\left\| u_{\lambda }\right\| _{L^{\infty }(\Bbb{T)}%
} $. Let $FC$ be the closed subspace of $FL^{\infty }$ consisting of those
sequences $\left\{ \lambda _{n}\right\} _{n\in \Bbb{Z}}$ which are the
Fourier coefficients of a continuous function on $\Bbb{T}$, i.e., a function
in $C(\Bbb{T)}$.

We will use the following basic properties of Fourier series
(\cite{Ka}, I.2 and I.3): 

\begin{fact}
  \begin{thmenumerate}
\label{cffs} 
\item
for each sequence $\lambda =$ $\left\{ \lambda _{n}\right\}
_{n\in \Bbb{Z}}\in FL^{1}$ the (Laurent or trigonometric) polynomials $%
P_{N}(e^{it})=\sum_{n=-N}^{N}\left( 1-\frac{\left| n\right| }{N+1}\right)
e^{int}\lambda _{n}$ converge to the associated function $u_{\lambda
}(e^{it})$ for a.e.\ $t\in [0,2\pi ]$, and also in the norm of $L^{1}(\Bbb{T)%
}$.

\item If, furthermore, $\lambda \in FL^{\infty }$, then $\sup_{t\in [0,2\pi
)}\left| P_{N}(e^{it})\right| \le \left\| u_{\lambda }\right\| _{L^{\infty }(%
\Bbb{T)}}$.

\item If $\lambda \in FC$, then the convergence of the above sequence $%
\left\{ P_{N}(e^{it})\right\} _{N\in \Bbb{N}}$ occurs for every $t\in
[0,2\pi ]$ and furthermore it is uniform on $\Bbb{T}$, i.e., the sequence
converges in $L^{\infty }(\Bbb{T)}$ norm.	
  \end{thmenumerate}
\end{fact}

For $p=1,\infty $ and each fixed $\alpha \in \Bbb{R}$ let $FL_{\alpha }^{p}$
be the space of complex sequences $\left\{ \lambda _{n}\right\} _{n\in \Bbb{Z%
}}$ such that $\left\{ e^{\alpha n}\lambda _{n}\right\} _{n\in \Bbb{Z}}\in
FL^{p}$. It is normed by $\left\| \left\{ \lambda _{n}\right\} _{n\in \Bbb{Z}%
}\right\| _{FL_{\alpha }^{p}}=\left\| \left\{ e^{\alpha n}\lambda
_{n}\right\} _{n\in \Bbb{Z}}\right\| _{FL^{p}}$. We also define $FC_{\alpha
} $ to be the closed subspace of $FL_{\alpha }^{\infty }$ of sequences $%
\left\{ \lambda _{n}\right\} _{n\in \Bbb{Z}}$ such that $\left\{ e^{\alpha
n}\lambda _{n}\right\} _{n\in \Bbb{Z}}\in FC$.

\smallskip 
These spaces are used extensively in \cite{Ja}, in particular to
define the complex interpolation method as both a minimal and as a maximal
method.

Obviously $FL^{\infty }\subset $ $FL^{1}\subset \ell ^{\infty }$ and both
these inclusions are continuous. So, for $j=0,1$ we have the continuous
inclusions $FL_{j}^{\infty }\subset FL_{j}^{1}\subset \ell _{\min
\{1,e^{n}\}}^{\infty }=\left\{ \left\{ \lambda _{n}\right\} _{n\in \Bbb{Z}%
}:\sup_{n\in \Bbb{Z}}\min \{1,e^{n}\}\left| \lambda _{n}\right| <\infty
\right\} $. Consequently, $\vec{FL^{p}}=(FL_{0}^{p},FL_{1}^{p})$ is a Banach
couple for $p=1$, and for $p=\infty $.

\smallskip

The continuous inclusion $FL_{\theta }^{1}\subset [\vec{FL^{1}}]_{\theta }$
was shown in \cite{Ja}, and it is not difficult to also obtain the reverse
continuous inclusion. As explained in \cite{CKM}, the formula 
\begin{equation}
\lbrack \vec{FL^{\infty }}]_{\theta }=FC_{\theta }\text{ to within
equivalence of norms}  \label{tweon}
\end{equation}
can be deduced from $[\vec{FL^{1}}]_{\theta }=FL_{\theta }^{1}$ with the
help of a duality argument. We need (\ref{tweon}) crucially in this paper,
and since existing proofs of it in the literature are not completely
explicit, we provide a new proof below, in Appendix  \ref{apzeon}. 
In fact there we prove the following slightly stronger result:

\begin{proposition}
\label{zeon}For each $\theta \in (0,1)$ the three spaces $[\vec{FL^{\infty }}%
]_{\theta }$, $[\vec{FL^{\infty }}]_{\theta ,\Bbb{A}}$ and $FC_{\theta }$
coincide isometrically.
\end{proposition}

Recall that the equality $[\vec A]_{\theta,\Bbb A} = [\vec A]_\theta$ 
holds for all Banach couples $\vec A$ \cite{Cw1}, but in general the norms are
equivalent only and not identical. (At least, this is proven for thick
annuli in \cite{Cw1}; it is undoubtedly true for all annuli, but it
seems that no-one has published a proof of this.)

\smallskip

We next consider some special and natural maps on the spaces $FL_{\alpha
}^{p}$ and $FC_{\alpha }$.

\begin{definition}
\label{amush}\smallskip For each fixed $s\in \Bbb{R}$, let $\mathcal{M}_{s}$
be the ``multiplier'' map on sequences defined by $\mathcal{M}_{s}\left\{
\lambda _{n}\right\} _{n\in \Bbb{Z}}=\left\{ e^{ins}\lambda _{n}\right\}
_{n\in \Bbb{Z}}$ and for each fixed $m\in \Bbb{Z}$ let $\mathcal{S}_{m}$ be
the ``shift'' map on sequences defined by $\mathcal{S}_{m}\left\{ \lambda
_{n}\right\} _{n\in \Bbb{Z}}=\left\{ \lambda _{n-m}\right\} _{n\in \Bbb{Z}}$.
\end{definition}

\begin{lemma}
\label{bmush}For each $s\in \Bbb{R}$ and each $m\in \Bbb{Z}$ both $\mathcal{M%
}_{s}$ and $\mathcal{S}_{m}$ are bounded maps of $FL_{\alpha }^{p}$ onto
itself and of $FC_{\alpha }$ onto itself, for each fixed $\alpha $ and for $%
p=1$, $\infty $. In fact $\mathcal{M}_{s}$ is an isometry, and $\mathcal{S}%
_{m}$ has norm $e^{\alpha m}$ on each of these spaces.
\end{lemma}

\begin{proof}
This is obvious from the following simple observation:
Whenever $\left\{ e^{\alpha n}\lambda _{n}\right\} _{n\in \Bbb{Z}}$ is the
sequence of Fourier coefficients of the function $v(e^{it})\in L^{1}(\Bbb{T)}
$, then

(i) $\{e^{ins}e^{\alpha n}\lambda _{n}\}_{n\in \Bbb{Z}}=\left\{ e^{\alpha
n}\left( \mathcal{M}_{s}\{\lambda _{n}\}\right) _{n}\right\} $ is the
sequence of Fourier coefficients of the ``rotated'' function $e^{it}\mapsto
v(e^{i(t+s)})$, and

(ii) $\{e^{\alpha (n-m)}\lambda _{n-m}\}_{n\in \Bbb{Z}}=e^{-\alpha m}\left\{
e^{\alpha n}\left( \mathcal{S}_{m}\{\lambda _{n}\}\right) _{n}\right\} $ is
the sequence of Fourier coefficients of the function $e^{it}\mapsto
e^{imt}v(e^{it})$. 
\end{proof}

The results in the following theorem and proposition are slight
reformulations of straightforward and classical properties of analytic
functions. But it seems preferable to state them explicitly in the notation
needed for our purposes, and to provide explicit proofs.

\begin{theorem}
\label{dtfyv}Suppose that the sequence $\lambda =\left\{ \lambda
_{n}\right\} _{n\in \Bbb{Z}}\in FL^{1}$ is also an element of $FL_{1}^{1}$.
Then $\lambda \in FC_{\theta }$ for all $\theta \in (0,1)$. Furthermore, the
limit 
\begin{equation}
v(z)=\lim_{N\rightarrow \infty }\sum_{n=-N}^{N}\left( 1-\frac{\left|
n\right| }{N+1}\right) z^{n}\lambda _{n}  \label{dvz}
\end{equation}
exists for every $z\in \Bbb{A}^{\circ }$ and almost every $z\in $ $\Bbb{T}$
and almost every $z\in e\Bbb{T}$, and the complex valued function $v$ which
it defines has the following properties:
\begin{romenumerate}
\item 
$v(e^{it})=u_{\lambda }(e^{it})$ for a.e.\ $t\in [0,2\pi )$,

\item 
the function $e^{it}\mapsto v(e^{1+it})$ is in $L^{1}(\Bbb{T)}$ and $%
\frac{1}{2\pi }\int_{0}^{2\pi }e^{-int}v(e^{1+it})dt=e^{n}\lambda _{n}$ for
all $n\in \Bbb{Z}$,

\item 
for each $\theta \in (0,1),$ the function $e^{it}\mapsto v(e^{\theta
+it})$ is in $C(\Bbb{T)}$ and \linebreakx
$\frac{1}{2\pi }\int_{0}^{2\pi
}e^{-int}v(e^{\theta +it})dt=e^{\theta n}\lambda _{n}$ for all $n\in \Bbb{Z}$%
, and

\item for each $\theta \in (0,1)$, 
\begin{equation}
\left| v(e^{\theta })\right| \le C_{\theta }\left( \int_{0}^{2\pi }\left|
v(e^{it})\right| dt\right) ^{1-\theta }\left( \int_{0}^{2\pi }\left|
v(e^{1+it})\right| dt\right) ^{\theta },  \label{urk}
\end{equation}
where the constant $C_{\theta }$ depends only on $\theta $.
\end{romenumerate}
\end{theorem}

\begin{proof} 
Since $\lambda \in FL^{1}\cap FL_{1}^{1}$, it follows that 
\begin{equation}
\left| \lambda _{m}\right| \le \frac{1}{2\pi }\min \{1,e^{-m}\}\max \left\{
\left\| \lambda \right\| _{FL^{1}(\Bbb{T)}},\left\| \left\{ e^{n}\lambda
_{n}\right\} _{n\in \Bbb{Z}}\right\| _{FL^{1}(\Bbb{T)}}\right\} \quad\text{for
each }m\in \Bbb{Z}.  \label{grerk}
\end{equation}
So, for each fixed $\theta \in (0,1)$, the sequence $\left\{ e^{(\theta
+it)n}\lambda _{n}\right\} _{n\in \Bbb{\Bbb{Z}}}$ satisfies 
\begin{equation}
\sup_{t\in \Bbb{R}}\left| e^{(\theta +it)n}\lambda _{n}\right| =\left|
e^{\theta n}\lambda _{n}\right| \le const.\min \{e^{\theta n},e^{(\theta
-1)n}\}\quad\text{for each }n\in \Bbb{Z}.  \label{tdp}
\end{equation}

For each $N\in \Bbb{N}$ we define the function (Laurent polynomial) $P_{N}:%
\Bbb{C}\backslash \{0\}\rightarrow \Bbb{C}$ by 
\begin{equation}
P_{N}(z)=\sum_{n=-N}^{N}\left( 1-\frac{\left| n\right| }{N+1}\right)
z^{n}\lambda _{n}.  \label{ylp}
\end{equation}

In view of Fact~\ref{cffs} (i), the sequence $\left\{ P_{N}(e^{it})\right\}
_{N\in \Bbb{N}}$ converges a.e.\ to $u_{\lambda }(e^{it})$ for a.e.\ $t\in
[0,2\pi )$ and also in the norm of $L^{1}(\Bbb{T)}$. For exactly the same
reason, the sequence $\left\{ e^{it}\mapsto P_{N}(e^{1+it})\right\} _{N\in 
\Bbb{N}}=\left\{ e^{it}\mapsto \sum_{n=-N}^{N}\left( 1-\frac{\left| n\right| 
}{N+1}\right) e^{(1+it)n}\lambda _{n}\right\} _{N\in \Bbb{N}}$ converges
a.e.\ on $\Bbb{T}$, and also in $L^{1}(\Bbb{T)}$ norm, to a function, (which
we shall temporarily call $u_{\lambda ,1}(e^{it})$) whose Fourier
coefficients are $\left\{ e^{n}\lambda _{n}\right\} $.

For each fixed $\theta \in (0,1)$ we see from (\ref{tdp}) that the sequence
of functions $e^{it}\mapsto $ $P_{N}(e^{\theta +it})$ converges absolutely
and uniformly on $\Bbb{T}$ to a continuous function on $\Bbb{T}$ (which we
shall temporarily call $u_{\lambda ,\theta }(e^{it})$) whose Fourier
coefficients are $\left\{ e^{\theta n}\lambda _{n}\right\} $.

\smallskip 
In view of the preceding remarks, it is clear that $\lambda \in
FC_{\theta }$ for each $\theta \in (0,1)$. 
It is also clear that the function 
$v$ defined as in (\ref{dvz}) i.e., by $v(z):=\lim_{N\rightarrow \infty
}P_{N}(z)$, exists for all $z\in \Bbb{A}^{\circ }$ and for almost every $z$
on each of the circles $\Bbb{T}$ and $e\Bbb{T}$, and it satisfies $%
v(e^{it})=u_{\lambda }(e^{it})$ and $v(e^{1+it})=u_{\lambda ,1}(e^{it})$ for
a.e.\ $t\in [0,2\pi )$ and $v(e^{\theta +it})=u_{\lambda ,\theta }(e^{it})$
for all $t\in [0,2\pi )$ and all $\theta \in (0,1)$. So parts (i), (ii) and
(iii) of Theorem \ref{dtfyv} have been established.

It still remains to prove part (iv). For this let us first remark that all
continuous functions $f:\Bbb{A}\rightarrow \Bbb{C}$ which are analytic in
the interior of $\Bbb{A}$ satisfy 
\begin{equation}
\left| f(e^{\theta })\right| \le C_{\theta }\left( \int_{0}^{2\pi }\left|
f(e^{it})\right| dt\right) ^{1-\theta }\left( \int_{0}^{2\pi }\left|
f(e^{1+it})\right| dt\right) ^{\theta }  \label{brqq}
\end{equation}
for each $\theta \in (0,1)$, where the constant $C_{\theta }$ depends only
on $\theta $. We can show this via classical results. Or we can deduce it as
a very very special case (where $X_{0}=X_{1}=\Bbb{C}$) of the estimate (1)
of part (i) of Lemma 2 of \cite{CwK} p.\ 265. (Note that in the proof of
part (i) of Lemma 2 of \cite[p. 265]{CwK}, there is a small misprint: the
function $F(z)=f(e^{z})e^{z^{2}}$ should of course be $%
F(z)=f(e^{z})e^{-z^{2}}$. Cf.\ also the proof of the related estimate $(**)$
on p.\ 363 of \cite{CKM}.)

Applying \smallskip (\ref{brqq}) to the analytic polynomials $P_{N}$, we
have 
\begin{equation}
\left| P_{N}(e^{\theta })\right| \le C_{\theta }\left( \int_{0}^{2\pi
}\left| P_{N}(e^{it})\right| dt\right) ^{1-\theta }\left( \int_{0}^{2\pi
}\left| P_{N}(e^{1+it})\right| dt\right) ^{\theta }.  \label{sayt}
\end{equation}

\smallskip Finally we use the fact that the above-mentioned convergence of $%
\left\{ P_{N}(e^{j+it})\right\} _{N\in \Bbb{N}}$ to $v(e^{j+it})$ in $L^{1}(%
\Bbb{T)}$ implies that $\lim_{N\rightarrow \infty }\int_{0}^{2\pi }\left|
P_{N}(e^{j+it})\right| dt=\int_{0}^{2\pi }\left| v(e^{j+it})\right| dt$ for $%
j=0,1$. Thus we obtain (\ref{urk}) from (\ref{sayt}) by taking the limit as $%
N$ tends to $\infty $. This completes the proof of Theorem \ref{dtfyv}. 
\end{proof}

\smallskip

{\tiny \smallskip }We conclude this section with a proposition which
paraphrases and expands upon some of the things in the statements and the
proofs of Theorem \ref{dtfyv} and Lemma \ref{bmush}. Since in our
application of these things in the next section we will only need to deal
with sequences $\lambda $ in the smaller space $FL^{\infty }\cap
FL_{1}^{\infty }$ instead of in $FL^{1}\cap FL_{1}^{1}$, we only consider
such sequences here. Also, it will be convenient for us to use a slightly
exotic variant (see (\ref{deft})) of the formula (\ref{dvz}) so that all the
functions that we have to deal with will be defined unambiguously at \textit{%
every} point of $\Bbb{A}$

\begin{proposition}
\smallskip \label{gpo}For each element $\lambda =\left\{ \lambda
_{n}\right\} _{n\in \Bbb{Z}}$ of the sequence space $FL^{\infty }\cap
FL_{1}^{\infty }$, let $T\lambda $ be the function defined by 
\begin{equation}
T\lambda (z)=
\begin{cases}
\lim_{N\rightarrow \infty }\sum_{n=-N}^{N}\left( 1-\frac{|n|}{N+1}\right)
z^{n}\lambda _{n} & \text{if this limit exists} \\ 
0 & \text{if the above limit does not exist}
\end{cases}
\label{deft}
\end{equation}
for each $z\in \Bbb{A}$.

Then:
\begin{romenumerate}
\item 
The formula (\ref{deft}) defines a map $T$ from $FL^{\infty }\cap
FL_{1}^{\infty }$ into a certain class of functions on $\Bbb{A}$, which is
``linear almost everywhere''. I.e., for each $\lambda $ and $\xi $ in $%
FL^{\infty }\cap FL_{1}^{\infty }$ and each $\alpha $ and $\beta $ in $\Bbb{C%
}$, 
\begin{equation}
T(\alpha \lambda +\beta \xi )(z)=\alpha T\lambda (z)+\beta T\xi (z)\text{
for all }z\in \Bbb{A}^{\circ }\text{ and for almost all }z\in \Bbb{T}\cup e%
\Bbb{T}.  \label{ael}
\end{equation}

Furthermore, for each $\lambda \in FL^{\infty }\cap FL_{1}^{\infty }$, the
function $v=T\lambda $ has the following properties:

\item  
For each $\alpha \in [0,1]$, $v$ satisfies 
\begin{equation}
\left\| \lambda \right\| _{FL_{\alpha }^{\infty }}=\underset{t\in [0,2\pi )}{%
\mathrm{ess~sup }}\left| v(e^{\alpha +it})\right| .  \label{aosf}
\end{equation}

\item  
When $\alpha =\theta \in (0,1)$ the function $v(e^{\theta +it})$
depends continuously on $t$ and \eqref{aosf} can be rewritten as 
\begin{equation}
\left\| \lambda \right\| _{FL_{\theta }^{\infty }}=\left\| \lambda \right\|
_{FC_{\theta }}=\sup_{t\in [0,2\pi )}\left| v(e^{\theta +it})\right| .
\label{isf}
\end{equation}

\item  
$v$ is analytic in $\Bbb{A}^{\circ }$.

\item  
The function $T(\mathcal{M}_{s}\lambda )$ coincides with the function $%
z\mapsto v(e^{is}z)$ for each fixed $s\in \Bbb{R}$.

\item  
Except possibly on some null subset of $\Bbb{T}\cup e\Bbb{T}$, the
function $T(\mathcal{S}_{m}\lambda )$ coincides with the function $z\mapsto
z^{m}v(z)$ for each fixed $m\in \Bbb{Z}$.

\item  
The function $v$ satisfies the estimate \eqref{urk}.
\end{romenumerate}
\end{proposition}

\smallskip

\smallskip 
\begin{proof} 
From Theorem \ref{dtfyv} and its proof it is
immediately clear that the right sides of (\ref{dvz}) and of (\ref{deft})
coincide for all $z\in \Bbb{A}^{\circ }$ and for almost all $z\in \Bbb{T}%
\cup e\Bbb{T}$ and that (keeping in mind the definitions of the spaces $%
FL_{\alpha }^{\infty }$ and $FC_{\alpha }$) properties (i), (ii), (iii) and
(vii) all hold.

To show property (iv) we note that, for any $\alpha $ and $\beta $
satisfying $0<\alpha <\beta <1$, and for all $z$ in the closed annulus $%
\left\{ z:e^{\alpha }\le \left| z\right| \le e^{\beta }\right\} $, the
sequence $\left\{ z^{n}\lambda _{n}\right\} _{n\in \Bbb{\Bbb{Z}}}$ satisfies
(cf.\ (\ref{grerk})) $\left| z^{n}\lambda _{n}\right| \le const.\min
\{1,e^{-n}\}\max \{e^{\alpha n},e^{\beta n}\}=const.\min \left\{ e^{\alpha
n},e^{(\beta -1)n}\right\} $. Consequently the sequence of Laurent
polynomials $\left\{ P_{N}(z)\right\} _{N\in \Bbb{N}}$ (cf.\ (\ref{ylp}))
converges uniformly on each such annulus and so the function $v=T\lambda
=\lim_{N\rightarrow \infty }P_{N}$ is indeed analytic in $\Bbb{A}^{\circ }$.

Property (v) is an immediate consequence of the definitions of $T$ and of $%
\mathcal{M}_{s}$ (cf.\ Definition \ref{amush}).

To show (vi), we choose an arbitrary $m\in \Bbb{Z}$ and let $w=T(\mathcal{S}%
_{m}\lambda )$. We have to show that 
\begin{equation}
w(z)=z^{m}v(z)\quad\text{for all }z\in \Bbb{A}^{\circ }\text{ and almost all }%
z\in \Bbb{T}\cup e\Bbb{T}.  \label{wzv}
\end{equation}
This is obvious for each $z\in \Bbb{A}^{\circ }$, since, when $1<\left|
z\right| <e$, we have 
\linebreakx
$\lim_{N\rightarrow \infty }\sum_{n=-N}^{N}\left( 1-%
\frac{|n|}{N+1}\right) z^{n}\lambda _{n}=\sum_{n=-\infty }^{\infty
}z^{n}\lambda _{n}$ and this series is absolutely convergent. It remains to
deal with $z\in \Bbb{T}\cup e\Bbb{T}$. When $z\in \Bbb{T}$, we apply part
(i) of Theorem \ref{dtfyv} to the sequence $\mathcal{S}_{m}\lambda $, and
see that the function $w=T(\mathcal{S}_{m}\lambda )$ must satisfy $\frac{1}{%
2\pi }\int_{0}^{2\pi }e^{-int}w(e^{it})dt=\lambda _{n-m}$ for each $n\in 
\Bbb{Z}$. Since we also have $\frac{1}{2\pi }\int_{0}^{2\pi
}e^{-int}[e^{imt}v(e^{it})]dt=\lambda _{n-m}$ for all $n\in \Bbb{Z}$, it
follows from the uniqueness theorem for Fourier series that 
\begin{equation}
w(e^{it})=e^{imt}v(e^{it})\quad\text{for a.e.\ }\ t\in [0,2\pi ).  \label{lcn}
\end{equation}
Similarly, from part (ii) of Theorem \ref{dtfyv} applied to both $\mathcal{S}%
_{m}\lambda $ and $\lambda $, we have 
\begin{equation*}
\frac{1}{2\pi }\int_{0}^{2\pi }e^{-int}w(e^{1+it})dt=e^{n}\lambda _{n-m}=%
\frac{1}{2\pi }\int_{0}^{2\pi }e^{-int}[e^{(1+it)m}v(e^{1+it})]dt
\end{equation*}
for all $n\in \Bbb{Z}$. So another application of the uniqueness theorem
gives us that 
\begin{equation}
w(e^{1+it})=e^{(1+it)m}v(e^{1+it})\quad\text{for a.e.\ }\ t\in [0,2\pi ).
\label{bcn}
\end{equation}
In view of (\ref{lcn}) and (\ref{bcn}) we have now established (\ref{wzv}).
This completes the proof of Proposition \ref{gpo}. 
\end{proof}

\smallskip

\section{The main result}

We are finally ready to state and prove our main result.

\begin{theorem}
\label{ikar}Let $\vec{B}=(B_{0},B_{1})$ be the couple $\vec{FL^{\infty }}$.
Let $K_{0}$ be a relatively compact subset of the unit ball of $FL^{\infty }$
and let $K_{1}$ be the unit ball of $\vec{FL_{1}^{\infty }}$. Let $\theta $
be a constant satisfying $0<\theta <1$ and let $K_{\theta }$ be the set of
all sequences $\lambda \in FL^{\infty }+FL_{1}^{\infty }$ of the form $%
\lambda =f(e^{\theta })$ for $f\in \mathcal{F}_{\Bbb{A}}(K_{0},K_{1})$. Then 
$K_{\theta }$ is a relatively compact subset of $FC_{\theta }$.
\end{theorem}

\begin{remark}
As stated in the introduction, Theorem \ref{ikar} immediately
applies an affirmative answer to Question C in the case where $(B_{0},B_{1})$
is the couple of sequence spaces $(FL^{\infty },FL_{1}^{\infty })$. This is
because of (\ref{tweon}) (cf.\ also Proposition \ref{zeon}). We will also
use (\ref{tweon}) in the proof of Theorem \ref{ikar}.
\end{remark}

\begin{proof} 
$K_{1}$ is of course closed and convex, and, since the
closed convex hull of a relatively compact subset of a normed linear space
is compact, we may also suppose without loss of generality that $K_{0}$ is
closed and convex. It then follows from Fact~\ref{konv} and (\ref{pgdensef})
that the set of sequences 
\begin{equation*}
K_{\theta ,\mathcal{G}}:=\{\lambda =f(e^{\theta }):f\in \mathcal{F}_{\Bbb{A}%
}(K_{0},K_{1})\cap \mathcal{G}_{\Bbb{A}}(\vec{FL^{\infty })\}}
\end{equation*}
is a dense subset of $K_{\theta }$ in the norm of $[\vec{FL^{\infty }}%
]_{\theta ,\Bbb{A}}=FC_{\theta }$. Thus we have reduced the proof of the
theorem to showing that the set $K_{\theta ,\mathcal{G}}$ is relatively
compact in $FC_{\theta }$.

Obviously $K_{\theta ,\mathcal{G}}\subset FL^{\infty }\cap FL_{1}^{\infty }$%
. Let $T$ be the map defined in (\ref{deft}) and let $T(K_{\theta ,\mathcal{G%
}})$ be the collection of all functions $v=T\lambda $ obtained as $\lambda $
ranges over $K_{\theta ,\mathcal{G}}$. Let $\left. T(K_{\theta ,\mathcal{G}%
})\right| _{e^{\theta }\Bbb{T}}$ denote the collection of all functions on
the circle $e^{\theta }\Bbb{T}$ which are restrictions to that circle of
functions in $T(K_{\theta ,\mathcal{G}})$. It follows from (\ref{isf}) and
the fact that $T\xi (z)-T\lambda (z)=T(\xi -\lambda )(z)$ for all $z\in
e^{\theta }\Bbb{T}$ and all $\xi ,\lambda \in FL^{\infty }\cap
FL_{1}^{\infty }$ (cf.\ (\ref{ael})) that, in order to show that $K_{\theta ,%
\mathcal{G}}$ is relatively compact in $FC_{\theta }$, it suffices to show
that $\left. T(K_{\theta ,\mathcal{G}})\right| _{e^{\theta }\Bbb{T}}$ is a
relatively compact subset of $C(e^{\theta }\Bbb{T)}$. By definition, $%
K_{\theta ,\mathcal{G}}$ is contained in the unit ball of $[\vec{FL^{\infty }%
}]_{\theta ,\Bbb{A}}$ and so $\left. T(K_{\theta ,\mathcal{G}})\right|
_{e^{\theta }\Bbb{T}}$ is a bounded set of functions in $C(e^{\theta }\Bbb{T)%
}$. So, by the Arzel\`{a}--Ascoli theorem, it remains only to show that $%
\left. T(K_{\theta ,\mathcal{G}})\right| _{e^{\theta }\Bbb{T}}$ is an
equicontinuous set of functions. Intuitively at least, this is a rather
obvious consequence of Theorem \ref{diag} and Theorem \ref{dtfyv}, in
particular (\ref{urk}). But let us carefully go through all the steps to
show it:

Consider some fixed but arbitrary element $\psi $ of $\left. T(K_{\theta ,%
\mathcal{G}})\right| _{e^{\theta }\Bbb{T}}$. It must be a continuous
function $\psi :e^{\theta }\Bbb{T}\rightarrow \Bbb{C}$ which is the
restriction to $e^{\theta }\Bbb{T}$ of some function $T\lambda $ where the
sequence $\lambda $ is a fixed element of $K_{\theta ,\mathcal{G}}$. 
It must, 
by the definition of $\mathcal{G}_{\Bbb{A}}$,
be of the form 
\begin{equation}
\lambda =\sum_{n=-M}^{M}e^{\theta n}\lambda ^{n}  \label{iuy}
\end{equation}
for some $M\in \Bbb{N}$, where $\lambda ^{n}\in FL^{\infty }\cap
FL_{1}^{\infty }$ for each $n$. (Here we are using a superscript index in $%
\lambda ^{n}$ to avoid confusion with the previous notation $\lambda _{n}$
for individual terms of the sequence $\lambda $.) Furthermore, we have $%
\sum_{n=-M}^{M}e^{(j+it)n}\lambda ^{n}\in K_{j}$ for $j=0,1$ and for each
fixed $t\in [0,2\pi )$. In fact, for these values of $j$ and $t$, we have 
\begin{equation}
\sum_{n=-M}^{M}e^{(j+it)n}\lambda ^{n}\in K_{j}\cap FL^{\infty }\cap
FL_{1}^{\infty }.  \label{bromb}
\end{equation}

Let $K_{*}=\left. T(K_{0}\cap FL^{\infty }\cap FL_{1}^{\infty })\right| _{%
\Bbb{T}}$, i.e., $K_{*}$ is the set of all 
functions defined on $\Bbb{T}$ which are the restrictions to $\Bbb{T}$
of functions of the form $T\xi $ for $\xi \in K_{0}\cap FL^{\infty }\cap
FL_{1}^{\infty }$. Let $K$ be the set of all elements (i.e., equivalence
classes of functions) in $L^{\infty }(\Bbb{T})$ which have a representative
in $K_{*}$. In view of (\ref{aosf}) for $\alpha =0$ and the fact that $T\xi
(z)-T\lambda (z)=T(\xi -\lambda )(z)$ for a.e.\ $z\in \Bbb{T} $ and all $\xi
,\lambda \in FL^{\infty }\cap FL_{1}^{\infty }$ (cf.\ (\ref{ael})), we can
assert that $K$ is a relatively compact subset of the unit ball of $%
L^{\infty }(\Bbb{T)}$. We now set $v_{n}=T(\lambda ^{n})$ for each $n\in
[-M,M]$, (where $\lambda ^{n}$ are the elements appearing in (\ref{iuy}))
and, for each $\zeta \in \Bbb{A}$, we define $u(\zeta )=\sum_{n=-M}^{M}\zeta
^{n}v_{n}$. Correspondingly, we define $u(\zeta ,w)=\sum_{n=-M}^{M}\zeta
^{n}v_{n}(w)$ for each $\zeta $ and $w$ in $\Bbb{A}$.

\smallskip

The following property of $u(\zeta ,w)$ is an obvious consequence of (\ref
{ael}): 
\begin{equation}
\text{For each fixed }\zeta \in \Bbb{A}\text{, }u(\zeta ,w)=T\left(
\sum_{n=-M}^{M}\zeta ^{n}\lambda ^{n}\right) (w)\text{ }
\begin{cases}
\text{(i)} & \text{for a.e.\ }w\in \Bbb{T}.   \\ 
\text{(ii)} & \text{for a.e.\ }w\in e\Bbb{T}.   \\ 
\text{(iii)} & \text{for all }w\in \Bbb{A}^{\circ }. 
\end{cases}
\label{wlcmp}
\end{equation}

Using part (i) of (\ref{wlcmp}) we see that, for each fixed $%
\zeta \in \Bbb{T}$, we can write $\zeta =e^{it}$ for some fixed real $t$ and
apply (\ref{bromb}) for the case where $j=0$ to deduce that the restriction
of the function $w\mapsto u(e^{it},w)$ to the circle $\Bbb{T}$ coincides for
a.e.\ $w\in \Bbb{T}$ with a function in $K_{*}$ (namely $T\left(
\sum_{n=-M}^{M}e^{itn}\lambda ^{n}\right) (w)$ restricted to $\Bbb{T}$). In
other words, if we let $\mathbf{v}_{n}$ denote the equivalence class in $%
L^{\infty }(\Bbb{T})$ containing $\left. v_{n}\right| _{\Bbb{T}}$ for each $%
n $, then, for each fixed $t\in [0,2\pi )$, the equivalence class $%
\sum_{n=-M}^{M}e^{int}\mathbf{v}_{n}$ in $L^{\infty }(\Bbb{T})$, which
contains the restriction of $u(e^{it},w)$ to $\Bbb{T}$, is an element of $K$.

This last condition means exactly that the restriction of $u(\zeta ,w)$ to $%
\Bbb{T}\times \Bbb{T}$ is an element of the class $\mathcal{E}_{K}$ defined
in Theorem \ref{diag}. Thus, if $\rho _{K}$ is the function defined by (\ref
{droq}), it follows from that definition that 
\begin{eqnarray}
&&\int_{0}^{2\pi }\left| u(e^{it},z_{1}e^{it})-u(e^{it},z_{2}e^{it})\right|
dt  \notag \\
&=&\int_{\Bbb{T}}|u(w,z_{1}w)-u(w,z_{2}w)|d\mu (w)\le \rho _{K}\left(
d(z_{1},z_{2})\right) \quad\text{for all }z_{1},z_{2}\in \Bbb{T}\text{.}
\label{vme}
\end{eqnarray}

Next, we consider the properties of the function $w\mapsto $ $u(\zeta ,w)$,
when $w$ is restricted to $e\Bbb{T}$ and for each constant $\zeta \in e\Bbb{T%
}$. This time we use part (ii) of (\ref{wlcmp}) to obtain that this function
coincides, for a.e.\ $w\in e\Bbb{T}$, with $T\left( \sum_{n=-M}^{M}\zeta
^{n}\lambda ^{n}\right) (w)$. By (\ref{bromb}) for $j=1$, we know that $%
\sum_{n=-M}^{M}\zeta ^{n}\lambda ^{n}$ is in the unit ball of $%
FL_{1}^{\infty }$. In turn this implies, using (\ref{aosf}) for $\alpha =1$,
that $\left| u(\zeta ,w)\right| \le 1$ for a.e.\ $w\in e\Bbb{T}$. 

Consider the function $U(\zeta ,w)$ defined for all 
$\zeta ,w\in \Bbb{T}$ by $U(\zeta ,w)=u(e\zeta ,ew)$. 
Clearly $U(\zeta ,w)=\sum_{n=-M}^{M}\zeta ^{n}V_{n}(w)$ where 
$V_{n}(w):=e^{n}v_{n}(ew)$.
We also have, for each fixed $%
\zeta \in \Bbb{T}$, that $\left| U(\zeta ,w)\right| \le 1$ for a.e.\ $w\in 
\Bbb{T}$. Thus Lemma~\ref{Lpddl} implies that for a.e.\ $w\in e\Bbb T$,
$|u(\zeta,w)|\le1$ for all $\zeta\in e\Bbb T$. It follows that for
every $z\in\Bbb T$, $|u(w,zw)|\le 1$ for a.e.\ $w\in e\Bbb T$, and
thus
\begin{equation}
\int_{0}^{2\pi }\left| u(e^{1+it},ze^{1+it})\right| dt\le 2\pi .
\label{tpz}
\end{equation}

The reason that we need (\ref{tpz}) is that it immediately implies that
\begin{equation}
\int_{0}^{2\pi }\left|
u(e^{1+it},z_{1}e^{1+it})-u(e^{1+it},z_{2}e^{1+it})\right| dt\le 4\pi, 
\quad\text{for all }z_{1},z_{2}\in \Bbb{T}.  \label{ree}
\end{equation}

\smallskip 

Now we will use part (iii) of (\ref{wlcmp}) with $\zeta
=e^{\theta }$. It gives us that 
\begin{equation}
\psi (w)=u(e^{\theta },w)=\sum_{n=-M}^{M}e^{\theta n}v_{n}(w)\quad\text{for all }%
w\in e^{\theta }\Bbb{T}.  \label{nsty}
\end{equation}

\smallskip 

Let $s$ and $\sigma $ be two arbitrary fixed real numbers and
consider the sequence $\xi =\sum_{n=-M}^{M}\mathcal{S}_{n}\left( \mathcal{M}%
_{s}-\mathcal{M}_{\sigma }\right) \lambda ^{n}$ where the $\lambda ^{n}$'s
are the same fixed sequences in $FL^{\infty }\cap FL_{1}^{\infty }$ which
were introduced in (\ref{iuy}). It follows from Lemma \ref{bmush} that $\xi
\in FL^{\infty }\cap FL_{1}^{\infty }$.

In view of parts (i), (v) and (vi) of Proposition \ref{gpo} and the
definitions of $v_{n}$ and $u(\zeta ,w)$ given above, we can assert that the
following sequence of equalities hold for all $w\in \Bbb{A}^{\circ }$ and
for almost all $w\in \Bbb{T}\cup e\Bbb{T}$. 
\begin{eqnarray}
T\xi (w) &=&\sum_{n=-M}^{M}T\left( \mathcal{S}_{n}\left( \mathcal{M}_{s}-%
\mathcal{M}_{\sigma }\right) \lambda ^{n}\right) (w)  \notag \\
&=&\sum_{n=-M}^{M}T\left( \mathcal{S}_{n}\mathcal{M}_{s}\lambda ^{n}\right)
(w)-\sum_{n=-M}^{M}T\left( \mathcal{S}_{n}\mathcal{M}_{\sigma }\lambda
^{n}\right) (w)  \notag \\
&=&\sum_{n=-M}^{M}w^{n}T\left( \mathcal{M}_{s}\lambda ^{n}\right)
(w)-\sum_{n=-M}^{M}w^{n}T\left( \mathcal{M}_{\sigma }\lambda ^{n}\right) (w)
\notag \\
&=&\sum_{n=-M}^{M}w^{n}T\lambda ^{n}(e^{is}w)-\sum_{n=-M}^{M}w^{n}T\lambda
^{n}(e^{i\sigma }w)  \notag \\
&=&\sum_{n=-M}^{M}w^{n}v_{n}(e^{is}w)-\sum_{n=-M}^{M}w^{n}v_{n}(e^{i\sigma
}w)  \notag \\
&=&u(w,e^{is}w)-u(w,e^{i\sigma }w).  \label{wsbjaf}
\end{eqnarray}

We now apply part (vii) of Proposition \ref{gpo} to the function $T\xi (w)$.
In view of (\ref{wsbjaf}), this gives us that 
\begin{multline*}
\left| u(e^{\theta },e^{is}e^{\theta })-u(e^{\theta },e^{i\sigma
}e^{\theta })\right|  \\
\le C_{\theta }\left( \int_{0}^{2\pi }\left|
u(e^{it},e^{is}e^{it})-u(e^{it},e^{i\sigma }e^{it})\right| dt\right)
^{1-\theta }
\cdot\\
\left( \int_{0}^{2\pi }\left|
u(e^{1+it},e^{is}e^{1+it})-u(e^{1+it},e^{i\sigma }e^{1+it})\right| dt\right)
^{\theta }.
\end{multline*}
In view of (\ref{nsty}), the first term in the preceding inequality is $%
\left| \psi (e^{\theta +is})-\psi (e^{\theta +i\sigma })\right| $. So, using
(\ref{vme}) and (\ref{ree}), we deduce that 
\begin{equation}
\left| \psi (e^{\theta +is})-\psi (e^{\theta +i\sigma })\right| \le
C_{\theta }\left( \rho _{K}(d(e^{is},e^{i\sigma })\right) ^{1-\theta }(4\pi
)^{\theta }\quad\text{for all }s,\sigma \in \Bbb{R}.  \label{urpp}
\end{equation}
The relative compactness of $K$ implies, via Theorem \ref{diag}, that $%
\lim_{\delta \rightarrow 0}\rho _{K}(\delta )=0$. Thus the inequality (\ref
{urpp}) establishes the equicontinuity of the set of functions $\left.
T(K_{\theta ,\mathcal{G}})\right| _{e^{\theta }\Bbb{T}}$. As explained
earlier, this suffices to complete the proof of Theorem \ref{ikar}. 
\end{proof}

\smallskip

\appendix
\section{\label{apzeon}A proof of Proposition \ref{zeon}.}

Until now we have not related very explicitly to Calder\'{o}n's
original definition of his spaces $[A_{0},A_{1}]_{\theta }$. We shall assume
that the reader is familiar with their construction via a certain space 
$\mathcal{F}(A_{0},A_{1})$ of functions $f:\Bbb{S}\rightarrow A_{0}+A_{1}$,
where $\Bbb{S}=\left\{ z\in \Bbb{C}:0\le \mathrm{Re}\,z\le 1\right\}$,
and with the convenient dense subspace 
$\mathcal{G}(A_{0},A_{1})$. We
refer to \cite{Ca} for the details. 

We will begin our proof of Proposition \ref{zeon}, by establishing the
inclusion 
\begin{equation}
FC_{\theta }\subset [FL^{\infty },FL_{1}^{\infty }]_{\theta }  \label{hpz}
\end{equation}
and the two norm inequalities 
\begin{equation}
\left\| \lambda \right\| _{[FL^{\infty },FL_{1}^{\infty }]_{\theta }}\le
\left\| \lambda \right\| _{FC_{\theta }}\quad\text{for all }\lambda \in
FC_{\theta },  \label{vya}
\end{equation}
and 
\begin{equation}
\left\| \lambda \right\| _{[FL^{\infty },FL_{1}^{\infty }]_{\theta ,\Bbb{A}%
}}\le \left\| \lambda \right\| _{FC_{\theta }}\quad\text{for all }\lambda \in
FC_{\theta }.  \label{sil}
\end{equation}

In fact it suffices to show that (\ref{vya}) and (\ref{sil}) each hold for
all finitely supported sequences $\lambda $. Since these form a dense subset
of $FC_{\theta }$, this immediately implies that they both hold for all $%
\lambda \in FC_{\Bbb{\theta }}$ and also establishes (\ref{hpz}).

We use essentially the same simple reasoning as was used in \cite{Ja} to
show that $FL_{\theta }^{1}\subset [FL^{1},FL_{1}^{1}]_{\theta }$. Given an
arbitrary sequence $\lambda =\left\{ \lambda _{n}\right\} _{n\in \Bbb{Z}}$
with finite support, we introduce the sequence valued function $f(z)=\left\{
e^{n(\theta -z)}\lambda _{n}\right\} _{n\in \Bbb{Z}}$. Clearly $f(\theta
)=\lambda $, and for all real $t$ and $j=0,1$ we have 
\begin{equation}\label{gzq}
  \begin{split}
\left\| f(j+it)\right\| _{FL_{j}^{\infty }} 
&=\left\| f(j+it)\right\|_{FC_{j}^{\infty }}
=\left\| \left\{ e^{nj+n(\theta -j-it)}\lambda_{n}
  \right\}_{n\in \Bbb{Z}}\right\|_{FC}
\\&
=\left\| \left\{ e^{n\theta
}e^{-int}\lambda _{n}\right\} _{n\in \Bbb{Z}}\right\| _{FC}  
=\left\| \left\{ e^{n\theta }\lambda _{n}
  \right\}_{n\in \Bbb{Z}}\right\|_{FC}
\\&
=\left\| \left\{ \lambda _{n}\right\} _{n\in \Bbb{Z}}\right\|_{FC_{\theta }}.  
  \end{split}
\end{equation}

Since $f(z)$ is an entire function taking values in a finite dimensional
subspace of $FL^{\infty }\cap FL_{1}^{\infty }$ and is bounded on $\Bbb{S}$,
it is clear that, for each $\delta >0$, the function $f_{\delta
}(z)=e^{\delta (z-\theta )^{2}}f(z)$ is an element of Calder\'{o}n's space $%
\mathcal{F}(FL^{\infty },FL_{1}^{\infty })$. So $\lambda =f_{\delta }(\theta
)\in [FL^{\infty },FL_{1}^{\infty }]_{\theta }$. 
By (\ref{gzq}) we have 
\begin{equation*}
  \begin{split}
\left\| \lambda \right\| _{\mathcal{[}FL^{\infty },FL_{1}^{\infty }]_{\theta
}}
&\le \left\| f_{\delta }\right\| _{\mathcal{F}(FL^{\infty },FL_{1}^{\infty
})}
\\&
=\sup \left\{ e^{\delta (j+it-\theta)^{2}}\|f(j+it)\|_{FL_j^\infty}
 :j=0,1;\,t\in\Bbb{R}\right\}
\\&
\le
\max \left\{ e^{\delta \theta ^{2}},e^{\delta (1-\theta )^{2}}\right\}
\left\| \left\{ \lambda _{n}\right\} _{n\in \Bbb{Z}}\right\| _{FC_{\theta
}}.	
  \end{split}
\end{equation*}
Since $\delta $ can be chosen arbitrarily small, this gives (\ref{vya}).

\smallskip 

Let us also consider the function $g:\Bbb{A}\rightarrow
FL^{\infty }\cap FL_{1}^{\infty }$ defined by $g(\zeta )=\left\{ e^{n\theta
}\zeta ^{-n}\lambda _{n}\right\} _{n\in \Bbb{Z}}$. Since $g(e^{z})=f(z)$ we
of course have $g(e^{\theta })=\lambda $ and also $g\in \mathcal{F}_{\Bbb{A}%
}(FL^{\infty },FL_{1}^{\infty })$ with norm 
\begin{equation*}
\left\| g\right\| _{\mathcal{F}_{\Bbb{A}}(FL^{\infty },FL_{1}^{\infty
})}=\sup \left\{ \left\| g(e^{j+it})\right\| _{FL_{j}^{\infty }}:j=0,1;\,t\in
[0,2\pi )\right\} =\left\| \left\{ \lambda _{n}\right\} _{n\in \Bbb{Z}%
}\right\| _{FC_{\theta }}.
\end{equation*}
This immediately gives us (\ref{sil}).

We now know that (\ref{vya}) and (\ref{sil}) hold for all finitely supported
sequences and so, as already explained, this establishes all of (\ref{hpz}),
(\ref{vya}) and (\ref{sil}).

\smallskip

Now we turn to proving the inclusion and inequalities which are the reverse
of (\ref{hpz}), (\ref{vya}) and (\ref{sil}) respectively. Again we will use
density properties, more explicitly, the facts that $\mathcal{G}(FL^{\infty
},FL_{1}^{\infty })$ is dense in $\mathcal{F}(FL^{\infty },FL_{1}^{\infty })$
and $\mathcal{G}_{\Bbb{A}}(FL^{\infty },FL_{1}^{\infty })$ is dense in $%
\mathcal{F}_{\Bbb{A}}(FL^{\infty },FL_{1}^{\infty })$. Our main step will be
to show that, for each $g\in \mathcal{G}(FL^{\infty },FL_{1}^{\infty })$, we
have $g(\theta )\in FC_{\theta }$ and 
\begin{equation}
\left\| g(\theta )\right\| _{FC_{\theta }}\le \left\| g\right\| _{\mathcal{F}%
(FL^{\infty },FL_{1}^{\infty })}.  \label{www}
\end{equation}
Analogously, for every $g\in \mathcal{G}_{\Bbb{A}}(FL^{\infty
},FL_{1}^{\infty })$ we will show that $g(e^{\theta })\in FC_{\theta }$ and 
\begin{equation}
\left\| g(e^{\theta })\right\| _{FC_{\theta }}\le \left\| g\right\| _{%
\mathcal{F}_{\Bbb{A}}(FL^{\infty },FL_{1}^{\infty })}.  \label{zwww}
\end{equation}

If $g\in \mathcal{G}(FL^{\infty },FL_{1}^{\infty })$ then $g(z)$ is a finite
sum $g(z)=\sum_{m=1}^{M}\phi _{m}(z)a^{m}$ where each $a^{m}=\left\{
a_{n}^{m}\right\} _{n\in \Bbb{Z}}$ is a sequence in the space $FL^{\infty
}\cap FL_{1}^{\infty }$ and each $\phi _{m}$ is a scalar valued entire
function which is bounded on $\Bbb{S}$. 
Applying Theorem \ref{dtfyv} to each sequence $a^{m}$, 
we see that $a^{m}\in FC_{\theta }$
for each $m$. 
(Note that $FL^\infty_j\subset FL^1_j$ and that the proof of Theorem
\ref{dtfyv} does not use Proposition \ref{zeon}. In fact, we only need
the simple argument at the beginning of the proof.) 
So obviously also $g(\theta )\in FC_{\theta }$.

\smallskip 

Given any sequence $\lambda =\left\{ \lambda _{n}\right\} _{n\in \Bbb{Z}}\in
FL^{\infty }$ and any other sequence $\sigma =\left\{ \sigma _{n}\right\}
_{n\in \Bbb{Z}}$ whose support $\Bbb{Z}_{\sigma }=\left\{ n\in \Bbb{Z}%
:\sigma _{n}\ne 0\right\} $ is a finite set, we of course have that the
finite sum $\sum_{n\in \Bbb{Z}}\lambda _{n}\sigma _{n}=\sum_{n\in \Bbb{Z}%
_{\sigma }}^{{}}\lambda _{n}\sigma _{n}=\frac{1}{2\pi }\int_{0}^{2\pi
}f_{\lambda }(e^{it})f_{\sigma }(e^{-it})dt$, where $f_{\lambda }$ and $%
f_{\sigma }$ are the functions whose sequences of Fourier coefficients are
respectively $\lambda $ and $\sigma $. In particular this means that 
\begin{equation}
\left| \sum_{n\in \Bbb{Z}}\lambda _{n}\sigma _{n}\right| \le \frac{1}{2\pi }%
\left\| \lambda \right\| _{FL^{\infty }}\left\| \sigma \right\| _{FL^{1}}
\label{xxp}
\end{equation}
and that 
\begin{equation}
\left\| \lambda \right\| _{FL^{\infty }}=\sup_{\sigma \in Q}\left|
\sum_{n\in \Bbb{Z}}\lambda _{n}\sigma _{n}\right|   \label{fya}
\end{equation}
where $Q$ is the set of all finitely supported sequences $\sigma $ such that 
$\left\| \sigma \right\| _{FL^{1}}=\int_{0}^{2\pi }\left| \sum_{n\in \Bbb{Z}%
_{\sigma }}\sigma _{n}e^{int}\right| dt=2\pi $.

\smallskip 

For each $\sigma \in Q$, we define the function $\psi _{\sigma }:%
\Bbb{C}\rightarrow \Bbb{C}$ by 
\begin{equation}
\psi _{\sigma }(z)=\sum_{m=1}^{M}\phi _{m}(z)\sum_{n\in \Bbb{Z}_{\sigma
}}a_{n}^{m}\sigma _{n}e^{nz}.  \label{oov}
\end{equation}
This is an entire function which is bounded on $\Bbb{S}$. So, by the
Phragm\`{e}n--Lindel\"{o}f theorem, we have that 
\begin{equation}
\left| \psi _{\sigma }(\theta )\right| \le \sup \left\{ \left| \psi _{\sigma
}(j+it)\right| :j=0,1,t\in \Bbb{R}\right\} .  \label{ppq}
\end{equation}
Applying (\ref{xxp}), we see that, for $j=0,1$ and for all $t\in \Bbb{R}$, 
\begin{align}
\left| \psi _{\sigma }(j+it)\right| 
&\le \frac{1}{2\pi }\left\| \left\{
\sum_{m=1}^{M}\phi _{m}(j+it)a_{n}^{m}e^{n(j+it)}\right\} _{n\in \Bbb{Z}%
}\right\| _{FL^{\infty }}\left\| \sigma \right\| _{FL^{1}}   
\notag\\&
=\left\| \left\{ \sum_{m=1}^{M}\phi _{m}(j+it)a_{n}^{m}e^{nj}\right\}
_{n\in \Bbb{Z}}\right\| _{FL^{\infty }}=\left\| \left\{ \sum_{m=1}^{M}\phi
_{m}(j+it)a_{n}^{m}\right\} _{n\in \Bbb{Z}}\right\| _{FL_{j}^{\infty }} 
\notag\\
&=\left\| g(j+it)\right\| _{FL_{j}^{\infty }}\le \left\| g\right\| _{%
\mathcal{F}(FL^{\infty },FL_{1}^{\infty })}.  \label{eej}
\end{align}

By (\ref{fya}) we have that 
\begin{equation}
  \begin{split}
\sup_{\sigma \in Q}\left| \psi _{\sigma }(\theta )\right| 
&=\left\| \left\{
\sum_{m=1}^{M}\phi _{m}(\theta )a_{n}^{m}e^{n\theta }\right\} _{n\in \Bbb{Z}%
}\right\| _{FL^{\infty }}
=\left\| \left\{ \sum_{m=1}^{M}\phi _{m}(\theta
)a_{n}^{m}\right\} _{n\in \Bbb{Z}}\right\| _{FL_{\theta }^{\infty }}
\\&
=\left\|
g(\theta )\right\| _{FL_{\theta }^{\infty }}.  	
  \end{split}
\label{jjo}
\end{equation}
Since $g(\theta )\in FC_{\theta }$ we have $\left\| g(\theta )\right\|
_{FL_{\theta }^{\infty }}=\left\| g(\theta )\right\| _{FC_{\theta }}$ and so
(\ref{www}) follows from (\ref{jjo}), (\ref{ppq}) and (\ref{eej}). The proof
of (\ref{zwww}) is almost the same. The only differences are that this time
the functions $\phi _{m}$ are each of the form $\phi _{m}(z)=z^{k_{m}}$ for
some $k_{m}\in \Bbb{Z}$, and we have to use the function $\psi _{\sigma
}(z)=\sum_{m=1}^{M}\phi _{m}(z)\sum_{n\in \Bbb{Z}_{\sigma }}a_{n}^{m}\sigma
_{n}z^{n}$. By the maximum modulus principle for functions which are
analytic on $\Bbb{A}$, we have $\left| \psi _{\sigma }(e^{\theta })\right|
\le \max \left\{ \left| \psi _{\sigma }(e^{j+it})\right| :j=0,1,t\in [0,2\pi
)\right\} $. Analogously to (\ref{eej}) this last expression is bounded from
above by $\left\| g\right\| _{\mathcal{F}_{\Bbb{A}}(FL^{\infty
},FL_{1}^{\infty })}$, and we obtain (\ref{zwww}).

\smallskip 

Finally, if $f$ is an arbitrary element of $\mathcal{F}%
(FL^{\infty },FL_{1}^{\infty })$ then we approximate it in $\mathcal{F}%
(FL^{\infty },FL_{1}^{\infty })$ norm by a sequence $\left\{ g_{n}\right\}
_{n\in \Bbb{N}}$ in $\mathcal{G}(FL^{\infty },FL_{1}^{\infty })$. Thus, by (%
\ref{www}), $\left\{ g_{n}(\theta )\right\} _{n\in \Bbb{N}}$ is a Cauchy
sequence in $FC_{\theta }$. Since it converges to $f(\theta )$ in the norm
of $FL^{\infty }+FL_{1}^{\infty }$ it follows that $f(\theta )\in FC_{\theta
}$ and $\left\| f(\theta )\right\| _{FC_{\theta }}\le \left\| f\right\| _{%
\mathcal{F}(FL^{\infty },FL_{1}^{\infty })}$. Consequently $[FL^{\infty
},FL_{1}^{\infty }]_{\theta }\subset FC_{\theta }$ and $\left\| \lambda
\right\| _{FC_{\theta }}\le \left\| \lambda \right\| _{[FL^{\infty
},FL_{1}^{\infty }]_{\theta }}$ for all $\lambda \in [FL^{\infty
},FL_{1}^{\infty }]_{\theta }$. Since $\mathcal{G}_{\Bbb{A}}(FL^{\infty
},FL_{1}^{\infty })$ is dense in $\mathcal{F}_{\Bbb{A}}(FL^{\infty
},FL_{1}^{\infty })$, an exactly analogous argument using (\ref{zwww}) shows
that $\left\| \lambda \right\| _{FC_{\theta }}\le \left\| \lambda \right\|
_{[FL^{\infty },FL_{1}^{\infty }]_{\theta ,\Bbb{A}}}$. 
(Alternatively, one could use the first case and the general fact 
that the inclusion
$[\vec A]_{\theta,\Bbb A} \to [\vec A]_\theta$ has norm 1, see \cite{Cw1}.)

This completes the proof of Proposition \ref{zeon}.

\smallskip \smallskip

\newcommand\jour{\emph}
\newcommand\vol{}

\end{document}